\documentclass[11pt,reqno]{article}
\usepackage[affil-it]{authblk}
\usepackage{amssymb}
\usepackage{amsmath}
\usepackage{amsthm}
\usepackage{amsfonts}
\usepackage{mathrsfs}
\usepackage{epsfig}
\usepackage{subfig}
\usepackage{graphics}
\usepackage{xcolor}
\usepackage{multirow}
\usepackage{diagbox}
\DeclareMathOperator*{\argmin}{arg\,min}
\newtheorem{theorem}{Theorem}[section]

\newtheorem{lemma}[theorem]{Lemma}
\usepackage{hyperref}
\providecommand{\keywords}[1]{\small \textbf{\textit{Keywords---}} #1}

\title{Graph Regularized Sparse $L_{2,1}$ Semi-Nonnegative Matrix Factorization for Data Reduction}

\author[1]{Anthony Rhodes}
\author[1,*]{Bin Jiang}
\author[2]{Jenny Jiang}

\affil[1]{Fariborz Maseeh Department of Mathematics and Statistics, Portland State University, Portland OR, USA}
\affil[2]{Goizueta Business School, Emory University, Atlanta GA, USA}
\affil[*]{Corresponding author: Bin Jiang, bjiang@pdx.edu}

\date{}

\begin{document}

\maketitle

\begin{abstract}
Non-negative Matrix Factorization (NMF) is an effective algorithm for multivariate data analysis, including applications to feature selection, pattern recognition, and computer vision. Its variant, Semi-Nonnegative Matrix Factorization (SNF), extends the ability of NMF to render parts-based data representations to include mixed-sign data. Graph Regularized SNF builds upon this paradigm by adding a graph regularization term to preserve the local geometrical structure of the data space. Despite their successes, SNF-related algorithms to date still suffer from instability caused by the Frobenius norm due to the effects of outliers and noise. In this paper, we present a new \textit{$L_{2,1}$} SNF algorithm that utilizes the noise-insensitive $L_{2,1}$ norm. We provide monotonic convergence analysis of the \textit{$L_{2,1}$} SNF algorithm. In addition, we conduct numerical experiments on three benchmark mixed-sign datasets as well as several randomized mixed-sign matrices to demonstrate the performance superiority of \textit{$L_{2,1}$} SNF over conventional SNF algorithms under the influence of Gaussian noise at different levels.
\end{abstract}


\maketitle

\section{Introduction}
\label{sec:introduction}

Data reduction algorithms represent an essential component of machine learning systems. The use of such reductions are supported by topological properties of data, including the well-known Manifold Hypothesis \cite{Narayanan2010}. Non-negative Matrix Factorization (NMF) \cite{Kong2011, Lee2001, Lin2007a} is a popular reduction algorithm that is defined as a problem of finding a matrix factorization of a given non-negative matrix $\mathbf{X}^{m\times n}$ so that $\mathbf{X}\approx \mathbf{UV}^T$ for non-negative factors $\mathbf{U}^{m \times k}$ and $\mathbf{V}^{n \times k}$. Data reduction can be achieved if $k$ is relatively small compared with $m$ and $n$ as the storage for $\mathbf{UV}^T$ tends to be smaller than that for $\mathbf{X}$, that is, $k(m+n) < mn$.

A significant feature of this methodology is that it gives rise to a parts-based decomposition of $\mathbf{X}$. Each column of $\mathbf{U}$ represents a basis element in the reduced space while the rows of $\mathbf{V}$ can be interpreted as embodying the coordinates for each basis element that render an approximation of the columns of $\mathbf{X}$. Since the number of basis vectors $k$ is often relatively small, this set of vectors represents a useful latent structure in the data. Finally, because each component in the factorization is restricted to be non-negative, their interaction in approximating $\mathbf{X}$ is strictly \textit{additive}, meaning that the columns of $\mathbf{U}$ yield a parts-based decomposition of $\mathbf{X}$.

There are several approaches to solve NMF problems. The seminal paper published in Nature by Lee and Seung \cite{Lee1999} in 1999 and their subsequent work \cite{Lee2001} provide a very effective multiplicative updates (MU) algorithm to solve NMF based on minimization of square of the following Euclidean distance (also regarded as Frobenius norm):
\begin{equation}
\|\mathbf{X}-\mathbf{UV}^T\|_{F}=\sqrt{\sum_{ij}\left(\mathbf{X}_{ij}-(\mathbf{UV}^T)_{ij}\right)^{2}}.
\label{l2norm}
\end{equation}
The MU algorithm becomes extremely popular due to its simple algorithm structure and has inspired a large amount of research work in the field of NMF and its application. To expedite the numerical performance of the MU algorithm, projected gradient method \cite{Lin2007b}, alternating nonnegative least squares \cite{Kim2008}, and alternating least squares \cite{Liu2013} are thus proposed for algorithm enhancement. Meanwhile, the MU algorithm has been extended with various new features and has also been used in solving real world application problems. Cai {\it et al.} \cite{Cai2011} propose a graph regularized NMF to encode the geometrical information of the data space with a nearest neighborhood graph for a new parts-based data representation preserving the geometrical structure of the data. Hoyer \cite{Hoyer2004} proposes a sparsely-constrained NMF by adding an extra $L_1$-norm penalty term for more sparse representation than standard NMF. As $L_1$-norm regularization is non-differentiable and not suitable for feature selection, Nie {\it et al.} \cite{Nie2010} and Yang {\it et al.} \cite{yang2011} use $L_{2,1}$-norm regularization terms for sparse feature selection. Wu {\it et al.} \cite{Wu2018} use $L_{2,1}$-norm and graph regularization together for a more robust NMF algorithm. Lee {\it et al.} \cite{Lee2010} and Xing {\it et al.} \cite{Xing2021} construct semi-supervised NMF to utilize partially labeled data to improve the clustering performance. Jiao {\it et al.} \cite{Jiao2020} and Yu {\it et al.} \cite{Yu2021} use hyper-graph regularized NMF for data process in the biomedical field.

When the data matrix $\mathbf{X}$ is not strictly non-negative, NMF is inapplicable. Nevertheless, in many common use cases, a parts-based decomposition is still a desideratum for data compression with mixed-sign data. Ding {\it et al.} \cite{Ding2010} introduce a useful compromise to this end with the Frobenius norm based Semi-Nonnegative Matrix Factorization (SNF), where $\mathbf{U}$ can have mixed-sign data as $\mathbf{X}$ while $\mathbf{V}$ is constrained to be non-negative to ensure the compressed decomposition of $\mathbf{X}$ still \textit{additive}. Luo and Peng \cite{Luo2017} introduce graph regularization and sparsity control terms into SNF for a more efficient algorithm. Shu {\it et al.} \cite{Shu2021} construct a deep SNF model with data elastic preserving (that is, retaining the intrinsic geometric structure information) by adding two graph regularizers in each layer. Jiang {\it et al.} \cite{Jiangw2020} develop a Robust SNF with adaptive graph regularization for gene representation. Rousset {\it et al.} \cite{Rousset2018} use SNF technique to generalize pattern recognition for single-pixel imaging. Peng {\it et al.} \cite{Peng2022} construct a new two-dimensional SNF algorithm to handle data clustering without converting 2D data to vectors in advance for a better data representation.

It is worthwhile to mention that solving $\displaystyle{\min_{\mathbf{U}\geq 0,\mathbf{V}\geq 0}}\|\mathbf{X}-\mathbf{UV}^T\|^2_{F}$ for a given nonnegative $\mathbf{X}$ has been shown as NP-hard \cite{Vavasis2010}, due to the non-negativity constraints on $\mathbf{U, V}$ and the fact that the objective function is non-convex for $(\mathbf{U}, \mathbf{V})$ altogether. Therefore, there exists no reported algorithms which can find an optimal solution within polynomial-time with respect to $m,n,k<\min(m,n)$ in general. However, the objective function $\|\mathbf{X}-\mathbf{UV}^T\|^2_{F}$ is convex for either $\mathbf{U}$ or $\mathbf{V}$ separately, which motivates researchers to construct iterative methods by finding an optimal $\mathbf{U}$ followed an optimal $\mathbf{V}$ consecutively in order to decrease the objective function at each iteration. Therefore, NMF methods do not lead to a global minimizer of the objective function theoretically. Similarly, solving $\displaystyle{\min_{\mathbf{U},\mathbf{V}\geq 0}}\|\mathbf{X}-\mathbf{UV}^T\|^2_{F}$ for SNF also turns out to be NP-hard \cite{Gillis2015}, as the nonnegative constraint still applies for $\mathbf{V}$. Therefore, SNF algorithms also try to minimize the objective function first for $\mathbf{U}$ then for $\mathbf{V}$ separately at each iteration to ensure the objective function decreases. A recent research monograph by Nicolas Gillis \cite{Gillis2020} provides a complete literature review of many popular NMF and SNF algorithms. It turns out that such approximate methods are very efficient in solving various kinds of application problems. It is obvious that the above objective functions have infinite number of global minimizers since ($\mathbf{UD}$, $\mathbf{VD}^{-1}$) achieves the same objective as ($\mathbf{U}$, $\mathbf{V}$) for any positive diagonal matrix $\mathbf{D}$. Whether NMF and SNF algorithms converge to one of those global optimizers still remains as an open question. Our paper relies on the existing theoretical framework to construct a noise-insensitive SNF algorithm which also demonstrates experimental convergence towards the optimal solution for several randomly generated matrices in practice.

It is well known that the Frobenius norm in \eqref{l2norm} is unstable with respect to noise and outliers \cite{Kong2011, Chachlakis2019}. Kong {\it et al.} \cite{Kong2011} and Jiang {\it et al.} \cite{Jiang2020} use $L_{2,1}$ norm to solve NMF, but their method cannot be extended to accommodate mixed-sign data. Chachlakis {\it et al.} \cite{Chachlakis2019} use $L_1$ norm for tensor decomposition with arbitrary data, but the usage of $L_1$ norm poses inherent limitations due to its non-differentiability.

In this paper, we present a new $L_{2,1}$ SNF algorithm for mixed-sign data. The algorithm is based on solving a constrained optimization problem for matrix reconstruction with respect to the following $L_{2,1}$ norm so as to control the growth of the outliers effectively.
\begin{equation}
\|\mathbf{X} -\mathbf{UV}^T\|_{2,1}=\sum_{j=1}^{n} \sqrt{\sum_{i=1}^{m}\left(\mathbf{X}_{ij}-(\mathbf{UV}^T)_{ij}\right)^{2}}.
\label{l21norm}
\end{equation}

To show the advantage of the $L_{2,1}$ norm \eqref{l21norm} over the Frobenius norm \eqref{l2norm}, we set $\mathbf{X}=[\mathbf{X}_1, \cdots \mathbf{X}_n]$, $\mathbf{V}=\left[\begin{array}{c}\mathbf{V}_1\\ \cdots \\ \mathbf{V}_n\end{array}\right]$ and let $\|\cdot\|_2$ be the Euclidean vector norm. Then
$$\|\mathbf{X}-\mathbf{UV}^T\|_{F}=\sqrt{\sum_{j=1}^{n} \|\mathbf{X}_j-\mathbf{U}\mathbf{V}^T_j\|^2_2} \; ; \qquad \|\mathbf{X} -\mathbf{UV}^T\|_{2,1}=\sum_{j=1}^{n} \|\mathbf{X}_j-\mathbf{U}\mathbf{V}^T_j\|_2 \; .$$
Since the error $\|\mathbf{X}_j-\mathbf{U}\mathbf{V}^T_j\|_2$ for each data point $\mathbf{X}_j$ in the Frobenius norm is squared, an outlier with a large error can easily dominate the objection function and lead to inaccurate approximation. On the contrary, such error for each data point in the $L_{2,1}$ norm is not squared, then the error caused by an outlier does not dominate the objective function and avoids an inaccurate approximation from the outlier. Therefore, $L_{2,1}$ norm is less prone to outliers than Frobenius norm.

We highlight four major benefits of the proposed algorithm as below:

(1) While the standard SNF measures the matrix decomposition error $\|\mathbf{X-UV}^T\|_F^2$ in the Euclidean space, $L_{2,1}$ SNF exploits the noise-insensitive property of $L_{2,1}$ norm for the new error $\|\mathbf{X-UV}^T\|_{2,1}$. Therefore, the algorithm is more stable and efficient than other SNF algorithms.

(2) While all conventional NMF and SNF algorithms with graph regularization use the Frobenius norm with fixed Laplacian matrix for all the iterations, $L_{2,1}$ SNF enforces updating all the non-zero entries of its Laplacian with latest decomposition information at each iteration for more accuracy while still keeping the same computational complexity. It is the first time to apply $L_{2,1}$ norm to graph regularization for matrix factorization.

(3) A monotonic convergence analysis is provided for $L_{2,1}$ SNF algorithm based on a \textit{proxy loss} function. It is the first time that a monotonic convergence analysis towards SNF algorithm under $L_{2,1}$ norm is reported in the literature. The result is consistent with its superior numerical performance on three benchmark mixed-sign datasets and several randomized matrices over conventional SNF algorithms.

(4) To address the theoretical concern about convergence of $L_{2,1}$ SNF algorithm, a numerical experiment is implemented to investigate whether the algorithm can indeed reach a global minimizer of the objective function. Two arbitrary matrices $\mathbf{U}$ and $\mathbf{V}\geq 0$ are generated randomly while $\mathbf{X}=\mathbf{UV}^T$ is enforced so that $\displaystyle{\min_{\mathbf{U}, \mathbf{V}\geq 0}}\|\mathbf{X} -\mathbf{UV}^T\|_{2,1}=0$. It is observed that for any randomized initial matrices $(\mathbf{U_0,V_0})$, $L_{2,1}$ SNF always converges to an optimal $(\it{U},\it{V})$ with $\|\mathbf{X} -\it{U}\it{V}^T\|_{2,1}\approx 0$. This experiment not only demonstrates practical advantage of $L_{2,1}$ SNF algorithm in finding a global minimizer of the objective function, but also motivates general theoretical investigation on whether the solutions obtained from all SNF algorithms are optimal solutions in minimizing $\mathbf{X}-\mathbf{UV}^T$, which is still an open question in numerical linear algebra.

The remainder of the paper is organized as follows: In Section \ref{sec:algorithm}, we give a brief review of some related algorithms, introduce the new $L_{2,1}$ SNF algorithm, and conduct complexity analysis over all these algorithms to show similar computational complexities among them. In Section \ref{sec:convergence}, we prove the monotonic convergence of the new algorithm. In Section \ref{sec:experiment}, we run numerical experiments on three benchmark mixed-sign datasets as well as several randomized mixed-sign matrices to demonstrate that $L_{2,1}$ SNF achieves better performance than other conventional SNF algorithms under the influence of Gaussian noise at different levels. Conclusion and future research consideration are provided in Section \ref{sec:conclusion}.

\section{Algorithm Construction and Complexity Analysis}
\label{sec:algorithm}

In this section, we first review some related algorithms, and then we introduce details of the $L_{2,1}$ SNF algorithm. Complexity analyses of all the algorithms are given to demonstrate similar complexity of $L_{2,1}$ SNF algorithm with other algorithms. To simplify our description, $\mathbf{X}_i$, $\mathbf{U}_i$ and  $\mathbf{V}_i$ stand for the \textit{i}-th columns of $\mathbf{X}$, $\mathbf{U}$ and the \textit{i}-th row of $\mathbf{V}$; $\mathbf{R}_+$ is the set of all non-negative real numbers; $\mathbf{D}$, $\mathbf{\hat{D}}$ and $\mathbf{\overline{D}}$ denote diagonal matrices; $tr(\mathbf{A})$ signifies the trace of matrix $\mathbf{A}$; $|\mathbf{A}|$ denotes the matrix with entries being the absolute values of the entries of $\mathbf{A}$.

It should be noted that among all the reported SNF algorithms towards mixed-sign data, the original SNF \cite{Ding2010} and GR SNF \cite{Luo2017} are regarded as conventional methods without relying on any specific property of the given dataset. The other SNF algorithms focus on either extending the method to handle multilevel or high dimensional data clustering, or solving real application problems arising from different fields. Therefore, for a fair comparison, we compare $L_{2,1}$ SNF with these two conventional algorithms to observe how $L_{2,1}$ SNF evolves from them for feature enhancement and performance improvement.

\subsection{NMF}
\label{subsec:nmf}

Given $\mathbf{X}\in \mathbf{R}^{m \times n}_+$ and $k<\min(m,n)$, NMF seeks $\mathbf{U}\in \mathbf{R}^{m \times k}_+$, whose columns represent the reduced basis vectors, and $\mathbf{V}\in \mathbf{R}^{n \times k}_+$, whose rows represent the coordinates of the new basis vectors for the columns of $\mathbf{X}$, such that $\mathbf{X}\approx \mathbf{UV}^T$. This can be formulated as a minimization of the objective function $\mathscr{J}(\mathbf{U},\mathbf{V})$:
\begin{equation}
\min_{\mathbf{U}\geq 0,\mathbf{V}\geq 0}\left( \|\mathbf{X}-\mathbf{UV}^T\|_F^2 = \mathscr{J}(\mathbf{U},\mathbf{V})\right).
\label{nmf}
\end{equation}

As solving \eqref{nmf} turns out to be NP-hard, Lee and Seung \cite{Lee2001} propose the following iterative algorithm to determine an optimal U followed by an optimal V consecutively at each iteration where $t=0,1,2,\cdots$ denotes the numbers of iteration steps:
\begin{align}
    \mathbf{U}_{ij}(t\!+\!1) & = \mathbf{U}_{ij}(t) \frac{\left(\mathbf{XV}(t)\right)_{ij}}{\left(\mathbf{U}(t)\mathbf{V}^T(t)\mathbf{V}(t)\right)_{ij}} \nonumber \\
     \mathbf{V}_{ij}(t\!+\!1) & = \mathbf{V}_{ij}(t) \frac{\left(\mathbf{X}^T\mathbf{U}(t\!+\!1)\right)_{ij}}{\left(\mathbf{V}(t)\mathbf{U}^T(t\!+\!1)\mathbf{U}(t\!+\!1)\right)_{ij}}.
\label{nmfsol}
\end{align}

They prove that the objective function in \eqref{nmf} decreases monotonically with the sequence $\{\mathbf{U(t), V(t)}\}$:
\begin{equation}
\mathscr{J}(\mathbf{U}(t\!+\!1),\mathbf{V}(t\!+\!1))\leq \mathscr{J}(\mathbf{U}(t\!+\!1),\mathbf{V}(t)) \leq \mathscr{J}(\mathbf{U}(t),\mathbf{V}(t)).
\label{objective}
\end{equation}

\subsection{SNF}
\label{subsec:snf}

Given $\mathbf{X}\in \mathbf{R}^{m \times n}$ and $k<\min(m,n)$, SNF seeks $\mathbf{U}\in \mathbf{R}^{m \times k}$ and $\mathbf{V}\in \mathbf{R}^{n \times k}_+$ to minimize the following objective function:
\begin{equation}
\min_{\mathbf{U},\mathbf{V}\geq 0} \left(\|\mathbf{X}-\mathbf{UV}^T\|_F^2=\mathscr{J}(\mathbf{U},\mathbf{V})\right).
\label{snf}
\end{equation}

Similarly, as solving \eqref{snf} is also NP-hard, Ding {\it et al.} \cite{Ding2010} propose the following iterative algorithm and show that the objective function in \eqref{snf} also decreases monotonically like \eqref{objective}:
\begin{align}
    \mathbf{U}(t\!+\!1) & = \mathbf{XV}(t)\left(\mathbf{V}^T(t)\mathbf{V}(t)\right)^{-1} \nonumber \\
     \mathbf{V}_{ij}(t\!+\!1) & = \mathbf{V}_{ij}(t) \sqrt{\frac{\left(\mathbf{X}^T\mathbf{U}(t\!+\!1)\right)^+_{ij}+\left[\mathbf{V}(t)\left(\mathbf{U}^T(t\!+\!1)\mathbf{U}(t\!+\!1)\right)^-\right]_{ij}}{\left(\mathbf{X}^T\mathbf{U}(t\!+\!1)\right)^{-}_{ij}+\left[\mathbf{V}(t)\left(\mathbf{U}^T(t\!+\!1)\mathbf{U}(t\!+\!1)\right)^+\right]_{ij}}}
\label{snfsol}
\end{align}
where the positive and negative parts of $\mathbf{A}$ are $\mathbf{A}^+_{ij}=\frac{|\mathbf{A}|_{ij}+\mathbf{A}_{ij}}{2}$, $\mathbf{A}^-_{ij}=\frac{|\mathbf{A}|_{ij}-\mathbf{A}_{ij}}{2}$.

It is clear that the minimization solutions of \eqref{nmf} and \eqref{snf} are not unique as ($\mathbf{UD}$, $\mathbf{VD^{-1}}$) achieves the same goal as ($\mathbf{U}$, $\mathbf{V}$) for any positive diagonal matrix $\mathbf{D}$. Therefore, it is anticipated that additional constraints are enforced for $\mathbf{U}$ and $\mathbf{V}$ in order to reach a unique optimal solution at each step.

\subsection{Graph Regularized SNF}
\label{subsec:snf-grs}

Given $\mathbf{X} \in \mathbf{R}^{m\times n}$ and $k<\min(m,n)$, GR SNF seeks $\mathbf{U} \in \mathbf{R}^{m\times k}$ and
$\mathbf{V} \in \mathbf{R}_+^{n\times k}$ to solve the following optimization problem:
\begin{equation}
\min_{\mathbf{U},\mathbf{V}\geq 0} \left(\|\mathbf{X}-\mathbf{UV}^T\|_F^2
+\alpha \sum_{i,j=1}^{n}\frac{1}{2}\|\mathbf{V_i-V_j}\|_2^2w_{ij} + \beta\|\mathbf{U}\|_{2,1}
=\mathscr{J}(\mathbf{U},\mathbf{V})\right)
\label{grsnf}
\end{equation}
where $\alpha \geq 0$ and $\beta \geq 0$ are tuning parameters. The second term encodes the local
geometrical structure of the data into the factorization through a nearest neighbor graph where $\mathscr{N}_i$ denotes the nearest $p$ neighboring column vectors of $\mathbf{X}$ for $\mathbf{X}_i$, and $w_{ij}=1$ if $j\in \mathscr{N}_i$ or $i\in \mathscr{N}_j$, and $w_{ij}=0$ otherwise. The third term adds $L_{2,1}$-norm constraint on the basis matrix $\mathbf{U}$ for sparse feature selection. It is trivial to show the second term can be represented as
\begin{equation}
\sum_{i,j=1}^{n}\frac{1}{2}\|\mathbf{V_i-V_j}\|_2^2w_{ij} = tr(\mathbf{V}^T\mathbf{LV}) = tr(\mathbf{V}^T\mathbf{\overline{D}V}) - tr(\mathbf{V}^T\mathbf{WV})
\label{grsnf2}
\end{equation}
where $\mathbf{L}=\mathbf{\overline{D}}-\mathbf{W}$ is the graph Laplacian with $\mathbf{W}=[w_{ij}]$, $\mathbf{\overline{D}}_{ii}=\displaystyle{\sum_{j\neq i}}w_{ij}$, $\mathbf{\overline{D}}_{ij} = 0 \enspace i\neq j$.

Luo and Peng \cite{Luo2017} propose the following iterative scheme to solve \eqref{grsnf}:
\begin{align}
    \mathbf{U}(t\!+\!1) & = \mathbf{XV}(t)\left(\beta\mathbf{\hat{D}(t)}+\mathbf{V}^T(t)\mathbf{V}(t)\right)^{-1} \nonumber \\
    \mathbf{V}_{ij}(t\!+\!1) & = \mathbf{V}_{ij}(t)\sqrt{\frac{\left(\mathbf{X}^T\mathbf{U}(t\!+\!1)\right)^+_{ij}+\left[\mathbf{V}(t)\left(\mathbf{U}^T(t\!+\!1)\mathbf{U}(t\!+\!1)\right)^-\right]_{ij}+\alpha\left(\mathbf{WV}(t)\right)_{ij}}
                             {\left(\mathbf{X}^T\mathbf{U}(t\!+\!1)\right)^-_{ij}+\left[\mathbf{V}(t)\left(\mathbf{U}^T(t\!+\!1)\mathbf{U}(t\!+\!1)\right)^+\right]_{ij} +\alpha\left(\mathbf{\overline{D}V}(t)\right)_{ij}}}
\label{grsnfsol}
\end{align}
where $\mathbf{\hat{D}_{ii}}(t) = \frac{0.5}{\|\mathbf{U}_i(t)\|_{2}}$.

GR SNF extends SNF by adding new graph regularization and sparsity enhancement features controlled by the parameters $\alpha$ and $\beta$. In fact, there is one more term $\|\mathbf{V}^T\mathbf{V}-\mathbf{I}\|_F^2$ used in \eqref{grsnf} for discriminative information constraint. As this term competes with the graph regularization term towards $\mathbf{V}$ to deteriorate the algorithm performance, we omit it here. The best approach to handle with discriminative information is the so-called semi-supervised algorithm \cite{Xing2021} with prior label information. Meanwhile, the convergence analysis of the scheme \eqref{grsnfsol} is not given in \cite{Luo2017}.

\subsection{$L_{2,1}$ SNF}
\label{subsec:l21snf}

Given $\mathbf{X} \in \mathbf{R}^{m\times n}$ and $k<\min(m,n)$, the proposed $L_{2,1}$ SNF algorithm seeks $\mathbf{U} \in \mathbf{R}^{m\times k}$ and
$\mathbf{V} \in \mathbf{R}_+^{n\times k}$ to solve the following optimization problem:
\begin{equation}
\min_{\mathbf{U},\mathbf{V}\geq 0} \left(\|\mathbf{X}-\mathbf{UV}^T\|_{2,1}
+\alpha \sum_{i<j} \|\mathbf{V}_i-\mathbf{V}_j\|_{2}w_{ij} + \beta\|\mathbf{U}\|_{2,1}
=\mathscr{J}(\mathbf{U},\mathbf{V})\right)
\label{l21snf}
\end{equation}
where $\alpha, \beta \geq 0$ are tuning parameters.

Since the first term $\|\mathbf{X}-\mathbf{UV}^T\|_{2,1}$ in \eqref{l21snf} is convex in either $\mathbf{U}$ or $\mathbf{V}$ only while not
convex in $(\mathbf{U}, \mathbf{V})$ altogether, it is unrealistic to construct an algorithm to find the global minimum of \eqref{l21snf}. Instead, we try
to minimize the objective function with one variable only while fixing the other one, iteratively. Therefore, $L_{2,1}$ SNF algorithm is
proposed for an iterative solution of \eqref{l21snf} as below:
\begin{align}
    \mathbf{U}(t\!+\!1) & = \mathbf{X}\mathbf{D}(t)\mathbf{V}(t)\left(\beta\mathbf{\hat{D}}(t)+\mathbf{V}^T(t)\mathbf{D}(t)\mathbf{V}(t)\right)^{-1} \nonumber \\
    \mathbf{V}_{ij}(t\!+\!1) & = \mathbf{V}_{ij}(t)\sqrt{\frac{\left[\mathbf{D}(t)\left(\mathbf{X}^T\mathbf{U}(t\!+\!1)\right)^+\right]_{ij}\!+\!\left[\mathbf{D}(t)\mathbf{V}(t)\left(\mathbf{U}^T(t\!+\!1)\mathbf{U}(t\!+\!1)\right)^-\right]_{ij}\!+\!\alpha\left(\mathbf{W}(t)\mathbf{V}(t)\right)_{ij}}
                             {\left[\mathbf{D}(t)\left(\mathbf{X}^T\mathbf{U}(t\!+\!1)\right)^-\right]_{ij}\!+\!\left[\mathbf{D}(t)\mathbf{V}(t)\left(\mathbf{U}^T(t\!+\!1)\mathbf{U}(t\!+\!1)\right)^+\right]_{ij}\!+\!\alpha\left(\mathbf{\overline{D}}(t)\mathbf{V}(t)\right)_{ij}}}
\label{l21snfsol}
\end{align}
where $\mathbf{D}(t)$, $\mathbf{\hat{D}}(t)$, $\mathbf{\overline{D}}(t)$ are diagonal matrices,
$\mathbf{D}_{ii}(t)=\frac{1}{\|\mathbf{X}_i-\mathbf{U}(t)\mathbf{V}_i^T(t)\|_2}$, $\mathbf{\hat{D}}_{ii}(t)=\frac{1}{\|\mathbf{U}_i(t)\|_2}$.
$\mathbf{L}=\mathbf{\overline{D}}(t)-\mathbf{W}(t)$ is a graph Laplacian with $\mathbf{W}_{ii}(t)=0$,
$\mathbf{W}_{ij}(t)=\frac{w_{ij}}{\|\mathbf{V}_i(t)-\mathbf{V}_j(t)\|_{2}} \enspace i\neq j$, $\mathbf{\overline{D}}_{ii}(t)=\displaystyle{\sum_{j\neq i}}\mathbf{W}_{ij}(t)$.

Compared with SNF and GR SNF algorithms, $L_{2,1}$ SNF fully utilizes the noise-insensitive $L_{2,1}$-norm towards the matrix decomposition term, the graph regularization term and the sparsity control term for a decomposition which is less sensitive to outliers. However, the adoption of $L_{2,1}$-norm makes it more challenging to determine the critical points for $\mathbf{U}$ and $\mathbf{V}$ towards the optimal solution compared with square of Frobenius norm which is well-suited for optimization. To resolve this dilemma, we introduce a new \textit{proxy loss} function based on three diagonal matrices $\mathbf{D}$, $\mathbf{\overline{D}}$ and $\mathbf{\hat{D}}$ for those three terms, respectively. The new function is similar to square of Frobenius norm so that a sequence of $\{\mathbf{U,V}\}$ can be constructed to optimize the new function at each step. Then we can show that the sequence also leads to the monotonic decrease of the objective function in \eqref{l21snf}.

To ensure algorithmic stability in calculating the matrices $\mathbf{D}$, $\mathbf{\hat{D}}$ and $\mathbf{W}$ when the denominator of some entry approaches to zero, we set a threshold $\epsilon=10^{-10}$. Whenever the denominator is less than $\epsilon$, it is reset to $\epsilon$ so as to avoid division by zero. In fact, such matrix formulation accelerates the clustering process. For example, during the $\mathbf{V}$ update stage, if $\|\mathbf{X}_i-\mathbf{U}\mathbf{V}_i^T\|_{2}$ becomes small enough for some $\mathbf{V}_i$, then the $i$-th row approximation performs better than other rows of $\mathbf{V}$. Therefore, a larger $\mathbf{D}_{ii}=\frac{1}{|\mathbf{X}_i-\mathbf{U}\mathbf{V}_i^T\|_{2}}$ value ensures update of $\mathbf{V}_i$ by \eqref{l21snfsol} relies more on itself than other rows.

\subsection{Computational Complexity Analysis}
\label{subsec:complexity}

In this subsection, we compare computational complexities of all the above algorithms where $\mathbf{X}\in \mathbf{R}^{m \times n}$, $\mathbf{U}\in \mathbf{R}^{m \times k}$, $\mathbf{V}\in \mathbf{R}^{n \times k}_+$, $k<\min(m,n)$, $p$ is the fixed number of nearest neighbors, and $t$ is the number of iterations.

Note that we need $mnk$ operations (addition and multiplication) to compute $\mathbf{XV}$. That is, the complexity for $\mathbf{XV}$ is $mnk$. Similarly, the complexity for $\mathbf{UV}^T\mathbf{V}$ is $nk^2+mk^2\leq 2mnk$ since $k<\min(m,n)$. Then the complexity to compute $\mathbf{U}$ as well as $\mathbf{V}$ in the updating rule \eqref{nmfsol} of NMF is $O(mnk)$. Therefore, the complexity of NMF algorithm is $O(tmnk)$.

Since the formula to compute $\mathbf{V}$ in the updating rule \eqref{snfsol} of SNF is almost the same as NMF except taking the positive and negative parts of the corresponding matrices, the complexity to compute $\mathbf{V}$ is still $O(mnk)$. However, the formula to compute $\mathbf{U}$ of SNF is different from that of NMF. Complexities of $\mathbf{XV}$, $\mathbf{V}^T\mathbf{V}$, $(\mathbf{V}^T\mathbf{V})^{-1}$ and multiplication of $\mathbf{XV}$ with $(\mathbf{V}^T\mathbf{V})^{-1}$ are $mnk$, $nk^2$, $k^3$ and $mk^2$, which are all bounded by $mnk$. Therefore, the complexity of SNF algorithm is $O(tmnk)$, same as that of NMF.

GR SNF has three more terms than SNF, the sparsity control term $\mathbf{\hat{D}}$ for $\mathbf{U}$, and the graph regularization terms $\mathbf{W}$ and $\mathbf{\overline{D}}$ for $\mathbf{V}$ in the updating rule \eqref{grsnfsol}. $\mathbf{\hat{D}}$ is a diagonal matrix with $k$ diagonal entries, each of which can be calculated with $m$ operations. Therefore, the complexity for $\mathbf{\hat{D}}$ is $O(mk)$. Meanwhile, the $p$-nearest neighborhood graph construction for $\mathbf{W}$ and $\mathbf{\overline{D}}$ needs one-time $O(n^2m)$ operations due to pairwise distance computation and sorting \cite{Cai2011}. Since $\mathbf{W}$ is a sparse matrix with $p$ nonzero elements on each row, we only need $O(npk)$ operations to compute $\mathbf{WV}$. Also since $\mathbf{\overline{D}}$ is a diagonal matrix, we need $O(nk)$ operations to compute $\mathbf{\overline{D}V}$. Therefore, the complexity of GR SNF algorithm is $O(tmnk+n^2m)$.

 Compared with GR SNF, $L_{2,1}$ SNF algorithm requires calculation of the diagonal matrix $\mathbf{D}$ and update of all the non-zero entries of $\mathbf{W}$ and $\mathbf{\overline{D}}$ during each iteration step in the updating rule \eqref{l21snfsol}. Firstly, $\mathbf{D}$ has $n$ diagonal entries, each of which can be calculated with $mk$ operations. Therefore, the construction of $\mathbf{D}$ has complexity $O(mnk)$. Meanwhile, computation of $\mathbf{DV}$ and $\mathbf{D}(\mathbf{X}^T\mathbf{U})$ only needs $O(nk)$ operations as $\mathbf{D}$ is diagonal. Secondly, $\mathbf{W}$ has $np$ non-zero entries, each of which can be updated with $k$ operations. Thus the complexity to update $\mathbf{W}$ is only $O(nkp)$. Similarly, update of $\mathbf{\overline{D}}$ needs $O(np)$ operations. Therefore, $L_{2,1}$ SNF algorithm has the same complexity as GR SNF, $O(tmnk+n^2m)$.

\section{Algorithm Convergence Analysis}
\label{sec:convergence}

In this section, we demonstrate the monotonic convergence of the objective $\mathscr{J}(\mathbf{U},\mathbf{V})$ in \eqref{l21snf} for the sequence $\{\mathbf{U(t), V(t)}\}$ constructed by $L_{2,1}$ SNF \eqref{l21snfsol}. In \ref{subsec:Uconvergence}, we show that updating $\mathbf{U}$ by formula \eqref{l21snfsol} with fixed $\mathbf{V}$ yields a decrease of \eqref{l21snf}, that is, $\mathscr{J}\left(\mathbf{U}(t\!+\!1),\mathbf{V}(t)\right) \leq \mathscr{J}\left(\mathbf{U}(t),\mathbf{V}(t)\right)$. In \ref{subsec:Vconvergence}, we show that updating $\mathbf{V}$ by \eqref{l21snfsol} with fixed $\mathbf{U}$ also yields a decrease of \eqref{l21snf}, that is, $\mathscr{J}\left(\mathbf{U}(t\!+\!1),\mathbf{V}(t+1)\right) \leq \mathscr{J}\left(\mathbf{U}(t+1),\mathbf{V}(t)\right)$. In \ref{subsec:UVconvergence}, we combine these two estimates to conclude that the objective sequence $\mathscr{J}(\mathbf{U}(t),\mathbf{V}(t))$ decreases and thus converges to its infimum $J_{\textbf{inf}}=\mathscr{J}\left(\it{U},\it{V}\right)$ for one pair $(\it{U},\it{V})$. In fact, it still remains as an open question whether $\mathscr{J}\left(\it{U}, \it{V}\right)$ reaches the global minimum of $\mathscr{J}(\mathbf{U},\mathbf{V})$ for all related NMF and SNF methods. In \ref{subsec:optimality}, we show that for the sequence $\{\mathbf{U(t), V(t)}\}$, the formula \ref{l21snfsol} is indeed optimal for each variable while the other is fixed during each iterative step, so as to reduce $\mathscr{J}(\mathbf{U},\mathbf{V})$ as much as possible.

Let us define a \textit{proxy loss} function as below:
\begin{equation}
\mathscr{L}(\mathbf{U},\mathbf{V})\!=\!tr[(\mathbf{X}\!-\!\mathbf{UV}^T)\mathbf{D}(\mathbf{X}\!-\!\mathbf{UV}^T)^{T}]
+\alpha tr[\mathbf{V}^T\mathbf{LV}]+\beta tr[\mathbf{U\hat{D}U}^T]
\label{proxyloss}
\end{equation}
where $\mathbf{D} \in \mathbb{R}^{n\times n}$, $\mathbf{D}_{ii} = \frac{1}{\|\mathbf{X}_i-\mathbf{U}\mathbf{V}_i^T\|_2}$,
$\mathbf{\hat{D}} \in \mathbf{R}^{k\times k}$,
$\mathbf{\hat{D}_{ii}} = \frac{1}{\|\mathbf{U}_i\|_2}$.
$\mathbf{L}=\mathbf{\overline{D}}\!-\!\mathbf{W}\!\in\! \mathbf{R}^{n\times n}$ is a graph Laplacian with
$$\mathbf{W}_{ij}=\left\{\begin{array}{cc}
  \frac{w_{ij}}{\|\mathbf{V}_i-\mathbf{V}_j\|_{2}} & i\neq j \\
  0 & i = j
\end{array},\right.
\hspace{6mm}
\mathbf{\overline{D}}_{ij}=\left\{\begin{array}{cc}
  \displaystyle{\sum_{l=1}^{n}}\mathbf{W}_{il} & i = j \\
  0 & i \neq j
\end{array}.\right.
$$

\subsection{Monotonic Decrease of $\mathscr{J}(\mathbf{U},\mathbf{V})$ for $\mathbf{U}$ update with fixed $\mathbf{V}$}
\label{subsec:Uconvergence}

First, we derive an iterative update formula for $\mathbf{U}$. The function in \eqref{proxyloss} can be rewritten as
\begin{equation}
\mathscr{L}(\mathbf{U},\mathbf{V}) = tr[\mathbf{XDX}^T] - 2tr[\mathbf{XDVU}^T]+tr[\mathbf{UV}^T\mathbf{DVU}^T]
 +\alpha tr[\mathbf{V}^T\mathbf{LV}]+\beta tr[\mathbf{U\hat{D}U}^T].
\label{proxyloss2}
\end{equation}

Then
\begin{equation*}
\nabla_{\mathbf{U}}\mathscr{L}=-2\mathbf{XDV} + 2\mathbf{UV}^T\mathbf{DV} + 2\beta\mathbf{U\hat{D}}.
\end{equation*}

Solving $\nabla_{\mathbf{U}}\mathscr{L}=0$ yields the update formula for $\mathbf{U}$:
\begin{equation}
\mathbf{U}=[\mathbf{XDV}][\beta\mathbf{\hat{D}}+\mathbf{V}^T\mathbf{DV}]^{-1}.
\label{wformula}
\end{equation}

We now consider \eqref{wformula} as an iterative update rule for $\mathbf{U}$ at step $t$, denoted by $\mathbf{U}(t+1)$, where $\mathbf{U}(t)$ and $\mathbf{V}(t)$ are regarded as fixed at the $t$-th update.
Then $\mathbf{D}(t)_{ii}=1/\|\mathbf{X}_i-\mathbf{U}(t)\mathbf{V}_i(t)^T\|_{2}$, $\mathbf{\hat{D}}(t)_{ii} = 1/\|\mathbf{U}_i(t)\|_{2}$, are also regarded as fixed.
Therefore, the iterative update for matrix $\mathbf{U}$ is given by:
\begin{equation}
\mathbf{U}(t+1)=[\mathbf{XD}(t)\mathbf{V}(t)][\beta \mathbf{\hat{D}}(t)+\mathbf{V}^T(t)\mathbf{D}(t)\mathbf{V}(t)]^{-1}.
\label{wtplus1}
\end{equation}

\begin{lemma}
\label{lem:Uupdate}
Let $\mathbf{U}(t)$ and $\mathbf{U}(t+1)$ represent consecutive updates for $\mathbf{U}$ as prescribed by \eqref{wtplus1}. Then $\mathscr{L}(\mathbf{U}(t+1),\mathbf{V}(t)) \leq \mathscr{L}(\mathbf{U}(t),\mathbf{V}(t))$, that is,
\begin{align*}
 & tr[(\mathbf{X}-\mathbf{U}(t+1)\mathbf{V}^T(t))\mathbf{D}(t)(\mathbf{X}-\mathbf{U}(t+1)\mathbf{V}^T(t))^{T}] +\alpha tr[\mathbf{V}^T(t)\mathbf{L}(t)\mathbf{V}(t)] \\
 & + \beta tr[\mathbf{U}(t+1)\mathbf{\hat{D}}(t)\mathbf{U}^T(t+1)] \\
 \leq & tr[(\mathbf{X}-\mathbf{U}(t)\mathbf{V}^T(t))\mathbf{D}(t)(\mathbf{X}-\mathbf{U}(t)\mathbf{V}^T(t))^{T})] +\alpha tr[\mathbf{V}^T(t)\mathbf{L}(t)\mathbf{V}(t)] \\
 &  + \beta tr[\mathbf{U}(t)\mathbf{\hat{D}}(t)\mathbf{U}^T(t)].
\end{align*}
\end{lemma}
\begin{proof}
Since $\mathbf{U}(t+1)$ is a stationary point of $\mathscr{L}(\mathbf{U},\mathbf{V})$ at step $t$, the optimality of the update formula \eqref{wtplus1} holds if we can demonstrate \eqref{proxyloss} is convex for $\mathbf{U}$.

To this end, $\frac{\partial\mathscr{L}}{\partial U_{ij}}$ is computed as:
\begin{align*}
\frac{\partial\mathscr{L}}{\partial U_{ij}} & = 2(\mathbf{UV}^T\mathbf{DV})_{ij}-2(\mathbf{XDV})_{ij}+2\beta(\mathbf{U\hat{D}})_{ij} \\
                                            & = 2\sum_{l=1}^{k} \mathbf{U}_{il}(\mathbf{V}^T\mathbf{DV})_{lj}-2(\mathbf{XDV})_{ij}+2\beta(\mathbf{U\hat{D}})_{ij}.
\end{align*}

The Hessian of $\mathscr{L}$ is consequently written as:
\begin{equation*}
\frac{\partial\mathscr{L}}{\partial U_{ij}\partial U_{pq}}
=2(\mathbf{V}^T\mathbf{DV}+\beta\mathbf{\hat{D}})_{qj}\delta_{ip} \quad
\quad 1\leq i,p \leq m \quad 1\leq j,q \leq k.
\end{equation*}

Therefore, the Hessian of $\mathscr{L}$ is a block diagonal matrix with each block in the form of $2\mathbf{V}^T\mathbf{DV}+2\beta \mathbf{\hat{D}}$.
It follows that $\mathscr{L}$ is convex for $\mathbf{U}$, thus the formula given for $\mathbf{U}$ in  \eqref{wtplus1} minimizes $\mathscr{L}$ in \eqref{proxyloss} at step $t$, i.e.,
$$\mathscr{L}(\mathbf{U}(t+1),\mathbf{V}(t)) \leq \mathscr{L}(\mathbf{U}(t),\mathbf{V}(t)).$$
\end{proof}

\begin{lemma}
\label{lem:inequality}
Given two arbitrary positive arrays $(x_1,x_2,\cdots, x_n)$ and $(y_1,y_2,\cdots, y_n)$, the following relation always holds:
$$\sum_{i=1}^n y_i - \sum_{i=1}^n x_i \leq \frac{1}{2}\left(\sum_{i=1}^n\frac{y_i^2}{x_i} - \sum_{i=1}^n \frac{x_i^2}{x_i}\right).$$
\end{lemma}
\begin{proof}
The above inequality is surely equivalent to the following:
$$\sum_{i=1}^n y_i \leq \frac{1}{2}\left(\sum_{i=1}^n\frac{y_i^2}{x_i} + \sum_{i=1}^n x_i\right).$$

For each $1\leq i \leq n$, we use the fact $\dfrac{a+b}{2}\geq \sqrt{ab}$ to obtain:
$$\frac{1}{2}\left(\frac{y_i^2}{x_i}+x_i\right)\geq \sqrt{\frac{y_i^2}{x_i}\cdot x_i}=y_i.$$

Summing them up yields the desired relation.
\end{proof}

\noindent \textbf{Remark.} Lemma~\ref{lem:inequality} plays an important role in proving monotonic decrease of the objective function in \eqref{l21snf} based on monotonic decrease
of $\mathscr{L}(\mathbf{U},\mathbf{V})$ in \eqref{proxyloss}, first for $\mathbf{U}$ update then for $\mathbf{V}$ update, as shown below. The left and right hand sides
of Lemma~\ref{lem:inequality} represent monotonicity of \eqref{l21snf} and \eqref{proxyloss}, respectively.

\begin{lemma}
\label{lem:Ucompare}
The following inequality holds:
\begin{align*}
 & \|\mathbf{X}-\mathbf{U}(t+1)\mathbf{V}^T(t)\|_{2,1}-\|\mathbf{X}-\mathbf{U}(t)\mathbf{V}^T(t)\|_{2,1}
  + \beta \|\mathbf{U}(t+1)\|_{2,1} - \beta \|\mathbf{U}(t)\|_{2,1} \\
 & \leq \frac{1}{2} tr\left[(\mathbf{X}-\mathbf{U}(t+1)\mathbf{V}^T(t))\mathbf{D}(t)(\mathbf{X}-\mathbf{U}(t+1)(\mathbf{V}^T(t))^T\right] \\
 & \hspace{4mm} - \frac{1}{2} tr\left[(\mathbf{X}-\mathbf{U}(t)\mathbf{V}^T(t))\mathbf{D}(t)(\mathbf{X}-\mathbf{U}(t)\mathbf{V}^T(t))^T\right] \\
 & \hspace{4mm} + \frac{\beta}{2} tr\left[\mathbf{U}(t+1)\mathbf{\hat{D}}(t)\mathbf{U}^T(t+1)\right]
                - \frac{\beta}{2} tr\left[\mathbf{U}(t)\mathbf{\hat{D}}(t)\mathbf{U}^T(t)\right].
\end{align*}
\end{lemma}
\begin{proof}
By setting $y_i=\|\mathbf{X}_i-\mathbf{U}(t+1)\mathbf{V}_i^T(t)\|_2$ and
$x_i=\|\mathbf{X}_i-\mathbf{U}(t)\mathbf{V}_i^T(t)\|_2$ in Lemma~\ref{lem:inequality},
we can see combination of the first two terms of left hand side $\leq$ combination of the first two terms of right hand side.

Similarly, by setting $y_i=\|\mathbf{U}_i(t+1)\|_2$ and $x_i=\|\mathbf{U}_i(t)\|_2$ in Lemma~\ref{lem:inequality},
we can see combination of the last two terms of left hand side $\leq$ combination of the last two terms of right hand side.

Combining them together finishes proof of Lemma~\ref{lem:Ucompare}.
\end{proof}

\begin{theorem}
\label{thm:Uupdate}
Updating $\mathbf{U}$ using formula \eqref{wtplus1} with fixed $\mathbf{V}$ yields a monotonic decrease in the objective function defined by \eqref{l21snf}, that is,
$\mathscr{J}\left(\mathbf{U}(t\!+\!1),\mathbf{V}(t)\right) \leq \mathscr{J}\left(\mathbf{U}(t),\mathbf{V}(t)\right)$.
\end{theorem}
\begin{proof}
By Lemma~\ref{lem:Uupdate}, the right hand side of Lemma~\ref{lem:Ucompare} $\leq 0$, so does the left hand side of Lemma~\ref{lem:Ucompare}, i.e.,
\begin{align*}
    & \|\mathbf{X}-\mathbf{U}(t+1)\mathbf{V}^T(t)\|_{2,1}
      + \beta \|\mathbf{U}(t+1)\|_{2,1} \\
\leq & \|\mathbf{X}-\mathbf{U}(t)\mathbf{V}^T(t)\|_{2,1}
      + \beta \|\mathbf{U}(t)\|_{2,1}.
\end{align*}

Meanwhile, the second term of \eqref{l21snf}, $\alpha \displaystyle{\sum_{i<j}} \|\mathbf{V}_i(t)-\mathbf{V}_j(t)\|_{2}w_{ij}$ is fixed during the $\mathbf{U}$ update. This finishes proof of Theorem~\ref{thm:Uupdate}.
\end{proof}

\subsection{Monotonic Decrease of $\mathscr{J}(\mathbf{U},\mathbf{V})$ for $\mathbf{V}$ update with fixed $\mathbf{U}$}
\label{subsec:Vconvergence}

Next we derive an iterative update formula for $\mathbf{V}\geq 0$; subsequently we prove monotonic convergence of this
update rule by showing that the \textit{proxy loss} $\mathscr{L}(\mathbf{U}, \mathbf{V})$ given in \eqref{proxyloss} decreases
monotonically for fixed $\mathbf{U}$. Since the last term of \eqref{proxyloss}, $\beta tr[\mathbf{U}^T\mathbf{\hat{D}U}]$, is
fixed during the $\mathbf{V}$ update, we ignore it here. Define the corresponding truncated \textit{proxy loss} function:
$$F(\mathbf{V})=tr[(\mathbf{X}-\mathbf{UV}^T)\mathbf{D}(\mathbf{X}-\mathbf{UV}^T)^{T}]
                +\alpha tr[\mathbf{V}^T\mathbf{LV}].$$

Based on $tr(\mathbf{AB})=tr(\mathbf{BA})$, $F(\mathbf{V})$ can be rewritten as:
\begin{equation}
F(\mathbf{V})\!=\!tr[\mathbf{XDX}^{T}]\!-\!2tr[\mathbf{V}^T\mathbf{DX}^T\mathbf{U}]\!+\!tr[\mathbf{U}^T\mathbf{UV}^T\mathbf{DV}]\!+\!\alpha tr[\mathbf{V}^T\mathbf{LV}].
\label{fhformula}
\end{equation}

To prove convergence for $\mathbf{V}$, we use an auxiliary function $\mathscr{A}(\mathbf{V},\mathbf{V'})$ as in \cite{Lee2001}.

$\mathscr{A}(\mathbf{V},\mathbf{V'})$ is defined as an auxiliary function for $F(\mathbf{V})$ if
\begin{equation}
F(\mathbf{V}) \leq \mathscr{A}(\mathbf{V},\mathbf{V'}), \quad F(\mathbf{V}) = \mathscr{A}(\mathbf{V},\mathbf{V}).
\label{auxiliary}
\end{equation}

\begin{lemma}
\label{lem:auxiliary}
If $\mathscr{A}(\mathbf{V},\mathbf{V'})$ is an auxiliary function of $F(\mathbf{V})$, then $F(\mathbf{V})$ is non-increasing under the update:
\begin{equation*}
\mathbf{V}^{t+1}=\argmin_{\mathbf{V}} \mathscr{A}(\mathbf{V}, \mathbf{V}^{t}).
\end{equation*}
\end{lemma}
\begin{proof}
\begin{equation*}
F(\mathbf{V}^{t+1})\leq
\mathscr{A}(\mathbf{V}^{t+1},\mathbf{V}^{t})\leq
\mathscr{A}(\mathbf{V}^{t},\mathbf{V}^{t})=F(\mathbf{V}^{t}).
\end{equation*}
\end{proof}
We now consider an explicit solution for $\mathbf{V}$ in the form of an iterative update, for which we subsequently prove convergence. Since $\mathbf{V}$ is non-negative, it is helpful to decompose both the $k \times k$ matrix $\mathbf{\Omega}=\mathbf{U}^T\mathbf{U}$ and the $n \times k$ matrix $\mathbf{\Phi}=\mathbf{X}^T\mathbf{U}$ into their positive and negative entries as below:
\begin{equation*}
\mathbf{\Omega}^{+}_{ij}=\frac{1}{2}(|\mathbf{\Omega}_{ij}|+\mathbf{\Omega}_{ij}), \quad
 \mathbf{\Omega}^{-}_{ij}=\frac{1}{2}(|\mathbf{\Omega}_{ij}|-\mathbf{\Omega}_{ij}).
\end{equation*}

\begin{lemma}
\label{lem:Vformula}
Under the iterative update:
\begin{equation}
\mathbf{V}_{ij}(t+1)\!=\!\mathbf{V}_{ij}(t)\sqrt{
\frac{(\mathbf{D}(t)\mathbf{\Phi}^+)_{ij}+(\mathbf{D}(t)\mathbf{V}(t)\mathbf{\Omega}^-)_{ij}
       +\alpha(\mathbf{W}(t)\mathbf{V}(t))_{ij}}
{(\mathbf{D}(t)\mathbf{\Phi}^-)_{ij}+(\mathbf{D}(t)\mathbf{V}(t)\mathbf{\Omega}^+)_{ij}
  +\alpha(\mathbf{\overline{D}}(t)\mathbf{V}(t))_{ij}}}
\label{Vformula}
\end{equation}
where $\mathbf{U}^{T}(t\!+\!1)\mathbf{U}(t\!+\!1)=\mathbf{\Omega} =\mathbf{\Omega}^{+}-\mathbf{\Omega}^{-}$, $\mathbf{X}^T\mathbf{U}(t\!+\!1)=\mathbf{\Phi}=\mathbf{\Phi}^{+}-\mathbf{\Phi}^{-}$,
the following relation holds for some auxiliary function $\mathscr{A}(\mathbf{V},\mathbf{V}')$ for $F(\mathbf{V})$:
\begin{equation*}
\mathbf{V}(t+1)=\argmin_{\mathbf{V}}
\mathscr{A}(\mathbf{V},\mathbf{V}(t)).
\end{equation*}
\end{lemma}
\begin{proof}
Using the notation introduced above, the truncated \textit{proxy loss} $F(\mathbf{V})$ in \eqref{fhformula} can be rewritten in the following form:
\begin{align}
F(\mathbf{V}) = \; & tr[\mathbf{XDX}^T]-2tr[\mathbf{V}^T\mathbf{D}\mathbf{\Phi}^+]+2tr[\mathbf{V}^T\mathbf{D}\mathbf{\Phi}^-] + tr[\mathbf{\Omega}^+\mathbf{V}^T\mathbf{DV}] \nonumber \\
                   & -tr[\mathbf{\Omega}^-\mathbf{V}^T\mathbf{DV}] + \alpha tr[\mathbf{V}^T\mathbf{\overline{D}V}]
                     - \alpha tr[\mathbf{V}^T\mathbf{WV}].
\label{fhterms}
\end{align}

In the subsequent steps, we provide an auxiliary function $\mathscr{A}(\mathbf{V},\mathbf{V}')$ for $F(\mathbf{V})$. Following \cite{Nie2010},
in order to construct an auxiliary function that furnishes an upper-bound for $F(\mathbf{V})$, we define $\mathscr{A}(\mathbf{V},\mathbf{V}')$ as a sum comprised of terms that represent upper-bounds for each of the positive terms and lower-bounds for each of the negative terms involving $\mathbf{V}$ in \eqref{fhterms}, respectively.

First, we derive a lower-bound for the second term of
\eqref{fhterms}, using $a \geq 1+log\; a, \; \forall \; a>0$:
\begin{equation}
tr[\mathbf{V}^T\mathbf{D}\mathbf{\Phi}^+]=\sum_{ij}\mathbf{V}_{ij} (\mathbf{D\Phi^+})_{ij} \geq
\sum_{ij}(\mathbf{D\Phi^+})_{ij}\mathbf{V}'_{ij}(1+\log\frac{\mathbf{V}_{ij}}{\mathbf{V}'_{ij}}).
\label{fh2ndterm}
\end{equation}

Second, using the fact that $a\leq \frac{a^2+b^2}{2b} \; \forall \; a,b>0$, we derive an upper bound for the third term of \eqref{fhterms}:
\begin{equation}
tr[\mathbf{V}^T\mathbf{D}\mathbf{\Phi}^-]=\sum_{ij}\mathbf{V}_{ij} (\mathbf{D\Phi^-})_{ij}\leq
 \sum_{ij}(\mathbf{D\Phi^-})_{ij}
 \frac{(\mathbf{V}_{ij})^2+(\mathbf{V'}_{ij})^2}{2\mathbf{V'}_{ij}}.
 \label{fh3rdterm}
\end{equation}

Third, we consider the fourth term of \eqref{fhterms}. Based on the facts that $\mathbf{D}$ is diagonal, $\mathbf{\Omega^+}$ is symmetric, and $2\sqrt{ab}\leq a+b$, we get:
\begin{align*}
tr[\mathbf{\Omega}^+\mathbf{V}^T\mathbf{DV}]
    & = tr[\mathbf{V}^T\mathbf{DV}\mathbf{\Omega}^+]
      =\sum_{ij}(\mathbf{DV}\mathbf{\Omega}^+)_{ij}\mathbf{V}_{ij} \\
    & = \sum_{ijk}(\mathbf{DV})_{ik}\mathbf{\Omega}^+_{kj}\mathbf{V}_{ij}
      = \sum_{ijk}\mathbf{D}_{ii}\mathbf{V}_{ik}\mathbf{\Omega}^+_{kj}\mathbf{V}_{ij}. \\
    & = \sum_{i,j<k}\mathbf{D}_{ii}\mathbf{\Omega}^+_{kj}2\mathbf{V}_{ik}\mathbf{V}_{ij}
         + \sum_{i,j=k}\mathbf{D}_{ii}\mathbf{\Omega}^+_{jj}\mathbf{V}_{ij}\mathbf{V}_{ij} \\
    & \leq \sum_{i,j<k}\mathbf{D}_{ii}\mathbf{\Omega}^+_{kj}\left(\frac{\mathbf{V'}_{ik}\mathbf{V}^2_{ij}}{\mathbf{V'}_{ij}}
         + \frac{\mathbf{V'}_{ij}\mathbf{V}^2_{ik}}{\mathbf{V'}_{ik}} \right)
         + \sum_{i,j=k}\mathbf{D}_{ii}\mathbf{\Omega}^+_{jj}\frac{\mathbf{V'}_{ij}\mathbf{V}^2_{ij}}{\mathbf{V'}_{ij}}.
\end{align*}

Since the second term inside the first summation corresponds to the case of $j>k$ due to symmetry of $\mathbf{\Omega^+}$, we then have:
\begin{equation}
tr[\mathbf{\Omega}^+\mathbf{V}^T\mathbf{DV}]
\leq \sum_{ijk}\frac{\mathbf{D}_{ii}\mathbf{V'}_{ik}\mathbf{\Omega}^+_{kj}\mathbf{V}^2_{ij}}{\mathbf{V'}_{ij}}
= \sum_{ij}\frac{(\mathbf{DV'}\mathbf{\Omega}^+)_{ij}\mathbf{V}^2_{ij}}{\mathbf{V'}_{ij}}.
\label{fh4thterm}
\end{equation}

Next, we consider the fifth term of \eqref{fhterms}:
\begin{align*}
tr[\mathbf{\Omega}^-\mathbf{V}^T\mathbf{DV}]
    & = tr[\mathbf{V}^T\mathbf{DV}\mathbf{\Omega}^-]
      =\sum_{ij}(\mathbf{DV}\mathbf{\Omega}^-)_{ij}\mathbf{V}_{ij} \\
    & = \sum_{ijk}(\mathbf{DV})_{ik}\mathbf{\Omega}^-_{kj}\mathbf{V}_{ij}
      = \sum_{ijk}\mathbf{D}_{ii}\mathbf{V}_{ik}\mathbf{\Omega}^-_{kj}\mathbf{V}_{ij}.
\end{align*}

We then use the inequality $a \geq 1+\log a$ again to obtain:
\begin{equation}
tr[\mathbf{\Omega}^-\mathbf{V}^T\mathbf{DV}]\geq
\sum_{ijk}\mathbf{D}_{ii}\mathbf{V'}_{ik}\mathbf{\Omega}^-_{kj}\mathbf{V'}_{ij}\left(1+\log\frac{\mathbf{V}_{ik}\mathbf{V}_{ij}}{\mathbf{V'}_{ik}\mathbf{V'}_{ij}}
\right).
\label{fh5thterm}
\end{equation}

Finally, we address the sixth and seventh terms of \eqref{fhterms} with the same techniques as for the fourth and fifth terms and obtain:
\begin{equation}
\alpha tr[\mathbf{V}^T\mathbf{\overline{D}V}] \leq \sum_{ij}\frac{\alpha(\mathbf{\overline{D}V'})_{ij}\mathbf{V}^2_{ij}}{\mathbf{V'}_{ij}}.
\label{fh6thterm}
\end{equation}
\begin{equation}
\alpha tr[\mathbf{V}^T\mathbf{WV}]\geq \sum_{ijk}\alpha \mathbf{W}_{ik}\mathbf{V'}_{kj}\mathbf{V'}_{ij}\left(1+\log\frac{\mathbf{V}_{kj}\mathbf{V}_{ij}}{\mathbf{V'}_{kj}\mathbf{V'}_{ij}}\right).
\label{fh7thterm}
\end{equation}

By setting all the formulas \eqref{fh2ndterm}, \eqref{fh3rdterm}, \eqref{fh4thterm}, \eqref{fh5thterm}, \eqref{fh6thterm}, and \eqref{fh7thterm} together,
we define the following auxiliary function $\mathscr{A}(\mathbf{V},\mathbf{V'})$:
\begin{align*}
\mathscr{A}(\mathbf{V},\mathbf{V'})
    & = tr[\mathbf{XD}\mathbf{X}^{T}]
      + 2\sum_{ij}(\mathbf{D\Phi^-})_{ij}\frac{(\mathbf{V}_{ij})^2+(\mathbf{V}'_{ij})^2}{2\mathbf{V'}_{ij}}
      + \sum_{ij}\frac{(\mathbf{DV'\Omega^+})_{ij}\mathbf{V}^2_{ij}}{\mathbf{V'}_{ij}} \\
    & -\!2\sum_{ij}(\mathbf{D\Phi^+})_{ij}\mathbf{V'}_{ij}\!\left(\!1\!+\!\log\frac{\mathbf{V}_{ij}}{\mathbf{V'}_{ij}}\!\right)\!
      \!-\! \sum_{ijk}\mathbf{D}_{ii}\mathbf{V'}_{ik}\mathbf{\Omega}^-_{kj}\mathbf{V'}_{ij}
        \!\left(\!1\!+\!\log\frac{\mathbf{V}_{ik}\mathbf{V}_{ij}}{\mathbf{V'}_{ik}\mathbf{V'}_{ij}} \!\right) \\
    & + \sum_{ij}\frac{\alpha (\mathbf{\overline{D}V'})_{ij}\mathbf{V}^2_{ij}}
                    {\mathbf{V'}_{ij}}
      - \sum_{ijk}\alpha \mathbf{W}_{ik}\mathbf{V'}_{kj}\mathbf{V'}_{ij}
                  \left(1+\log\frac{\mathbf{V}_{kj}\mathbf{V}_{ij}}{\mathbf{V'}_{kj}\mathbf{V'}_{ij}} \right).
\end{align*}

Observe that $F(\mathbf{V}) \leq \mathscr{A}(\mathbf{V},\mathbf{V'})$ and $F(\mathbf{V})=\mathscr{A}(\mathbf{V},\mathbf{V})$, as required for an auxiliary function in \eqref{auxiliary},
where $F(\mathbf{V})$ denotes the truncated \textit{proxy loss} as defined in \eqref{fhterms}. Lemma~\ref{lem:auxiliary} implies that $F(\mathbf{V})$ is non-increasing under the update:
$\mathbf{V}(t+1)=\argmin_{\mathbf{V}} \mathscr{A}(\mathbf{V},\mathbf{V}(t))$.

We now demonstrate that the minimum of $\mathscr{A}(\mathbf{V},\mathbf{V'})$ coincides with the update rule in \eqref{Vformula}. We first determine the critical points of $\mathscr{A}(\mathbf{V},\mathbf{V'})$:
\begin{align}
\frac{\partial\mathscr{A}(\mathbf{V},\mathbf{V'})}{\partial\mathbf{V}_{ij}}
    & = 2(\mathbf{D\Phi^-})_{ij}\Big(\frac{\mathbf{V}_{ij}}{\mathbf{V}'_{ij}}\Big)
      + 2\frac{(\mathbf{DV'\Omega^+})_{ij}\mathbf{V}_{ij}}{\mathbf{V'}_{ij}}
      + 2\alpha (\mathbf{\overline{D}V'})_{ij}\Big(\frac{\mathbf{V}_{ij}}{\mathbf{V}'_{ij}}\Big) \nonumber \\
    & -\!2(\mathbf{D\Phi^+})_{ij}\left(\frac{\mathbf{V'}_{ij}}{\mathbf{V}_{ij}}\right)
      \!-\!2\frac{(\mathbf{DV'\Omega^-})_{ij}\mathbf{V'}_{ij}}{\mathbf{V}_{ij}}
      \!-\!2\alpha (\mathbf{WV'})_{ij}\left(\frac{\mathbf{V'}_{ij}}{\mathbf{V}_{ij}}\right) \nonumber \\
    & = 0.
\label{ahhprimederiv}
\end{align}

Solving \eqref{ahhprimederiv} for $\mathbf{V}_{ij}$ yields the update formula given in \eqref{Vformula}. Thus \eqref{Vformula} corresponds with a critical point for $\mathscr{A}(\mathbf{V},\mathbf{V'})$.
Furthermore, as the Hessian of $\mathscr{A}(\mathbf{V},\mathbf{V'})$,
$$
\frac{\partial\mathscr{A}(\mathbf{V},\mathbf{V'})}{\partial \mathbf{V}_{ij}\partial \mathbf{V}_{kl}}\!=\!
\left\{
\begin{array}{ll}
     2\frac{(\mathbf{D\Phi^-})_{ij}}{\mathbf{V'}_{ij}}
      +2\frac{(\mathbf{DV'\Omega^+})_{ij}}{\mathbf{V}'_{ij}}
      +2\alpha \frac{(\mathbf{\overline{D}V'\Omega^+})_{ij}}{\mathbf{V}'_{ij}}      & \\
     +2\frac{(\mathbf{D\Phi^+})_{ij}\mathbf{V'}_{ij}}{\mathbf{V}_{ij}^2}
      \!+\!2\frac{(\mathbf{DV'\Omega^-})_{ij}\mathbf{V'_{ij}}}{\mathbf{V}_{ij}^2}
      \!+2\!\alpha \frac{(\mathbf{WV'})_{ij}\mathbf{V'_{ij}}}{\mathbf{V}_{ij}^2}    & \text{if } (i,j)\!==\!(k,l) \\
     0                                                                              & \text{otherwise}
\end{array}
\right.
$$
is positive definite, $\mathscr{A}(\mathbf{V},\mathbf{V'})$ is then convex for $\mathbf{V}$. Therefore, the critical point for $\mathscr{A}(\mathbf{V},\mathbf{V'})$ by \eqref{Vformula}
achieves the minimum of $\mathscr{A}(\mathbf{V},\mathbf{V'})$. This concludes proof of Lemma~\ref{lem:Vformula}.
\end{proof}

\noindent \textbf{Remark.} $\mathbf{V}\geq 0$ and $\mathbf{V'}\geq 0$ serve as the old and new matrices for the iterative update in Lemma~\ref{lem:Vformula}. \eqref{Vformula} indicates $\mathbf{V'}_{ij}$ is a multiple of $\mathbf{V}_{ij}$ component-wisely. If some $\mathbf{V}_{ij}=0$, then $\mathbf{V'}_{ij}=0$ as well. Therefore, the corresponding terms tend to zero on both sides of \eqref{fh2ndterm}, \eqref{fh3rdterm}, \eqref{fh4thterm}, \eqref{fh5thterm}, \eqref{fh6thterm}, and \eqref{fh7thterm} which hold as before.

\begin{lemma}
\label{lem:Vupdate}
Let $\mathbf{V}(t)$ and $\mathbf{V}(t+1)$ represent consecutive updates for $\mathbf{V}$ as prescribed by \eqref{Vformula}. Under this updating rule, the following inequality holds:
\begin{align*}
 & tr\left[(\mathbf{X}\!-\!\mathbf{U}(t\!+\!1)\mathbf{V}^T(t\!+\!1))\mathbf{D}(t)(\mathbf{X}\!-\!\mathbf{U}(t\!+\!1)\mathbf{V}^T(t\!+\!1))^{T}\right] \!+\! \alpha tr\left[\mathbf{V}(t\!+\!1)^T\mathbf{L}(t)\mathbf{V}(t\!+\!1)\right] \\
 & \leq tr\left[(\mathbf{X}\!-\!\mathbf{U}(t\!+\!1)\mathbf{V}^T(t))\mathbf{D}(t)(\mathbf{X}\!-\!\mathbf{U}(t\!+\!1)\mathbf{V}^T(t))^{T})\right] \!+\! \alpha tr\left[\mathbf{V}(t)^T\mathbf{L}(t)\mathbf{V}(t)\right].
\end{align*}
\end{lemma}
\begin{proof}
The proof follows directly from Lemma~\ref{lem:auxiliary} and Lemma~\ref{lem:Vformula}.
\end{proof}

\begin{lemma}
\label{lem:Vcompare}
The following inequality holds:
\begin{align*}
& \left(\|\mathbf{X}\!-\!\mathbf{U}(t\!+\!1)\mathbf{V}^T(t+1)\|_{2,1}
        \!+\!\alpha \sum_{i<j} \|\mathbf{V}_i(t+1)-\mathbf{V}_j(t+1)\|_{2}w_{ij}\right) \\
& - \left(\|\mathbf{X}\!-\!\mathbf{U}(t\!+\!1)\mathbf{V}^T(t)\|_{2,1}
         +\alpha \sum_{i<j} \|\mathbf{V}_i(t)-\mathbf{V}_j(t) \|_{2}w_{ij}\right) \\
& \leq \frac{1}{2}\left(tr\left[(\mathbf{X}\!-\!\mathbf{U}(t\!+\!1)\mathbf{V}^T(t\!+\!1))\mathbf{D}(t)(\mathbf{X}\!-\!\mathbf{U}(t)(\mathbf{V}^T(t\!+\!1))^T\right] \!+\! \alpha tr\left[\mathbf{V}(t\!+\!1)^T\mathbf{L}(t)\mathbf{V}(t\!+\!1)\right] \right) \\
& - \frac{1}{2}\left(tr\left[(\mathbf{X}\!-\!\mathbf{U}(t\!+\!1)\mathbf{V}^T(t))\mathbf{D}(t)(\mathbf{X}\!-\!\mathbf{U}(t)\mathbf{V}^T(t))^T\right] \!+\! \alpha tr\left[\mathbf{V}(t)^T\mathbf{L}(t)\mathbf{V}(t)\right] \right).
\end{align*}
\end{lemma}
\begin{proof}
By setting $y_i=\|\mathbf{X}_i-\mathbf{U}(t\!+\!1)\mathbf{V}_i^T(t+1)\|_2$ and $x_i=\|\mathbf{X}_i-\mathbf{U}(t\!+\!1)\mathbf{V}_i^T(t)\|_2$ in Lemma~\ref{lem:inequality},
we can see combination of the first and third terms of left hand side $\leq$ combination of the first and third terms of right hand side.

By setting $y_{ij}\!=\!\|\mathbf{V}_i(t\!+\!1)\!-\!\mathbf{V}_j(t\!+\!1)\|_{2}w_{ij}$ and $x_{ij}\!=\!\|\mathbf{V}_i(t)\!-\!\mathbf{V}_j(t)\|_{2}w_{ij}$ in Lemma~\ref{lem:inequality} again, and noticing
$$\sum_{i<j}\frac{y_{ij}^2}{x_{ij}} \!=\! \sum_{i<j}\frac{\|\mathbf{V}_i(t\!+\!1)\!-\!\mathbf{V}_j(t\!+\!1)\|_{2}^2\cdot w_{ij}}{\|\mathbf{V}_i(t)\!-\!\mathbf{V}_j(t)\|_{2}} \!=\! tr\left[\mathbf{V}(t\!+\!1)^T\mathbf{L}(t)\mathbf{V}(t\!+\!1)\right],$$

we can see combination of the second and fourth terms of left hand side $\leq$ combination of the second and fourth terms of right hand side.

Combining the above together finishes proof of Lemma~\ref{lem:Vcompare}.
\end{proof}

\begin{theorem}
\label{thm:Vupdate}
Updating $\mathbf{V}$ using formula \eqref{Vformula} with fixed $\mathbf{U}$ yields a monotonic decrease of the objective function in \eqref{l21snf}, that is,
$\mathscr{J}\left(\mathbf{U}(t\!+\!1),\mathbf{V}(t\!+\!1)\right) \leq \mathscr{J}\left(\mathbf{U}(t\!+\!1),\mathbf{V}(t)\right)$.
\end{theorem}
\begin{proof}
By Lemma~\ref{lem:Vupdate}, the right hand side of Lemma~\ref{lem:Vcompare} $\leq 0$, so does the left hand side of Lemma~\ref{lem:Vcompare}, i.e.,
\begin{align*}
&\|\mathbf{X}-\mathbf{U}(t\!+\!1)\mathbf{V}^T(t+1)\|_{2,1}
        +\alpha \sum_{i<j} \|\mathbf{V}_i(t+1)-\mathbf{V}_j(t+1)\|_{2}w_{ij} \\
& \leq \|\mathbf{X}-\mathbf{U}(t\!+\!1)\mathbf{V}^T(t)\|_{2,1}
         +\alpha \sum_{i<j} \|\mathbf{V}_i(t)-\mathbf{V}_j(t) \|_{2}w_{ij}.
\end{align*}

Meanwhile, the third term of \eqref{l21snf}, $\beta\|\mathbf{U}(t\!+\!1)\|_{2,1}$ is fixed during the $\mathbf{V}$ update. The proof is finished.
\end{proof}

\subsection{Monotonic Decrease of $\mathscr{J}(\mathbf{U},\mathbf{V})$ for $L_{2,1}$ SNF Algorithm}
\label{subsec:UVconvergence}
\begin{theorem}
\label{thm:UVupdate}
Let $(\mathbf{U}(t), \mathbf{V}(t))$ denote the $t$-th iterative value of $(\mathbf{U}, \mathbf{V})$ by $L_{2,1}$ SNF algorithm \eqref{l21snfsol}, then the objective sequence $\mathscr{J}(\mathbf{U}(t),\mathbf{V}(t))$ in \eqref{l21snf} decreases monotonically:
$$\mathscr{J}\left(\mathbf{U}(t\!+\!1),\mathbf{V}(t\!+\!1)\right) \leq \mathscr{J}\left(\mathbf{U}(t),\mathbf{V}(t)\right).$$
Therefore, it converges to its infimum $J_{\textbf{inf}}=\mathscr{J}\left(\mathbf{\it U},\mathbf{\it V}\right)$ for one pair $(\mathbf{\it U},\mathbf{\it V})$.
\end{theorem}
\begin{proof}
Theorem~\ref{thm:Uupdate} implies that
$$\mathscr{J}\left(\mathbf{U}(t\!+\!1),\mathbf{V}(t)\right) \leq \mathscr{J}\left(\mathbf{U}(t),\mathbf{V}(t)\right).$$

Meanwhile, Theorem~\ref{thm:Vupdate} implies that
$$\mathscr{J}\left(\mathbf{U}(t\!+\!1),\mathbf{V}(t\!+\!1)\right) \leq \mathscr{J}\left(\mathbf{U}(t\!+\!1),\mathbf{V}(t)\right).$$

Combination of the above two yields immediately
$$\mathscr{J}\left(\mathbf{U}(t\!+\!1),\mathbf{V}(t\!+\!1)\right) \leq \mathscr{J}\left(\mathbf{U}(t),\mathbf{V}(t)\right).$$
\end{proof}
Therefore, we prove the monotonic convergence of $L_{2,1}$ SNF algorithm by showing that $\mathscr{J}(\mathbf{U},\mathbf{V})$ decreases monotonically with the help of a smooth \textit{proxy loss} function $\mathscr{L}(\mathbf{U},\mathbf{V})$. The following Theorem reveals the relationship between the decreases of $\mathscr{J}(\mathbf{U},\mathbf{V})$ and $\mathscr{L}(\mathbf{U},\mathbf{V})$.
\begin{theorem}
\label{thm:JLdecrease}
\begin{equation}
\mathscr{J}\left(\mathbf{U}(t\!+\!1),\mathbf{V}(t)\right) - \mathscr{J}\left(\mathbf{U}(t),\mathbf{V}(t)\right) \leq \frac{1}{2}\left[\mathscr{L}\left(\mathbf{U}(t\!+\!1),\mathbf{V}(t)\right)-\mathscr{L}\left(\mathbf{U}(t),\mathbf{V}(t)\right) \right]
\label{Ucompare2}
\end{equation}
\begin{equation}
\mathscr{J}\left(\mathbf{U}(t\!+\!1),\mathbf{V}(t\!+\!1)\right) - \mathscr{J}\left(\mathbf{U}(t\!+\!1),\mathbf{V}(t)\right) \leq \frac{1}{2}\left[\mathscr{L}\left(\mathbf{U}(t\!+\!1),\mathbf{V}(t\!+\!1)\right)-\mathscr{L}\left(\mathbf{U}(t\!+\!1),\mathbf{V}(t)\right) \right].
\label{Vcompare2}
\end{equation}
\begin{equation}
\mathscr{J}\left(\mathbf{U}(t\!+\!1),\mathbf{V}(t\!+\!1)\right) - \mathscr{J}\left(\mathbf{U}(t),\mathbf{V}(t)\right) \leq \frac{1}{2}\left[\mathscr{L}\left(\mathbf{U}(t\!+\!1),\mathbf{V}(t\!+\!1)\right)-\mathscr{L}\left(\mathbf{U}(t),\mathbf{V}(t)\right) \right].
\label{UVcompare2}
\end{equation}
\end{theorem}
\begin{proof}
\eqref{Ucompare2} follows from Lemma~\ref{lem:Ucompare} while \eqref{Vcompare2} follows from Lemma~\ref{lem:Vcompare} directly.
Combination of \eqref{Ucompare2} and \eqref{Vcompare2} yields \eqref{UVcompare2} immediately.
\end{proof}
Therefore, the decrease of $\mathscr{J}(\mathbf{U},\mathbf{V})$ is at least half of the decrease of $\mathscr{L}(\mathbf{U},\mathbf{V})$ at each iteration. This is a reasonable estimate as $\mathscr{L}(\mathbf{U},\mathbf{V})$, which has three terms similar as square of Frobenius norm, tends to be larger than $\mathscr{J}(\mathbf{U},\mathbf{V})$, which has three terms with non-squared $L_{2,1}$ norm. Therefore, $L_{2,1}$ SNF algorithm provides an efficient way to reduce the objective function $\mathscr{J}(\mathbf{U},\mathbf{V})$ in \eqref{l21snf} at each iteration.

\subsection{Optimality of $L_{2,1}$ SNF Formula based on Convexity and KKT Condition}
\label{subsec:optimality}

Theorem~\ref{thm:UVupdate} shows that the objective sequence $\mathscr{J}(\mathbf{U}(t),\mathbf{V}(t))$ decreases and thus converges to its infimum $J_{\textbf{inf}}=\mathscr{J}\left(\it{U},\it{V}\right)$ for one pair $(\it{U},\it{V})$. However, it is still an open question that whether $\mathscr{J}\left(\it{U}, \it{V}\right)$ reaches the global minimum of $\mathscr{J}(\mathbf{U},\mathbf{V})$ for all related NMF and SNF methods \cite{Lee2001, Ding2010, Luo2017}. We now demonstrate step-wise optimality of our algorithm, that is, for the sequence $\{\mathbf{U(t), V(t)}\}$ constructed by \eqref{l21snfsol}, the objective $\mathscr{J}(\mathbf{U},\mathbf{V})$ in \eqref{l21snf} reaches its optimal value for each variable while the other is fixed during each iterative step.

Since the adoption of $L_{2,1}$-norm in the objective $\mathscr{J}(\mathbf{U},\mathbf{V})$ makes it difficult to determine the critical points for $\mathbf{U}$ and $\mathbf{V}$ at each step, we first show local optimality of our algorithm towards the \textit{proxy loss} function $\mathscr{L}(\mathbf{U}(t),\mathbf{V}(t)$ in \eqref{proxyloss}. Then we show that the \textit{proxy loss} $\mathscr{L}(\mathbf{U}(t),\mathbf{V}(t)$ and the objective $\mathscr{J}(\mathbf{U}(t),\mathbf{V}(t))$ go towards each other when $t\rightarrow \infty$, that is, $\mathscr{L}(\it{U},\it{V})=\mathscr{J}(\it{U},\it{V})$.
\begin{theorem}
\label{thm:Uoptimality}
Let $(\mathbf{U}(t), \mathbf{V}(t))$ denote the $t$-th iterative value of $(\mathbf{U}, \mathbf{V})$ by $L_{2,1}$ SNF algorithm \eqref{l21snfsol}, then
$$\mathscr{L}(\mathbf{U}(t+1),\mathbf{V}(t)) \leq \mathscr{L}(\mathbf{U},\mathbf{V}(t)) \qquad \forall \quad U\in \mathbf{R}^{m \times k}.$$
That is, $\mathscr{L}(\mathbf{U}(t+1),\mathbf{V}(t))=\displaystyle{\min_{\mathbf{U}}}\mathscr{L}(\mathbf{U},\mathbf{V}(t))$.
\end{theorem}
\begin{proof}
In Lemma~\ref{lem:Uupdate}, it has been shown that for fixed $\mathbf{V}(t)$, $\mathbf{U}(t+1)$ is a stationary point of $\mathscr{L}(\mathbf{U},\mathbf{V}(t))$ with respect to $\mathbf{U}$ at step $t$. Meanwhile, $\mathscr{L}(\mathbf{U},\mathbf{V}(t))$ is convex for $\mathbf{U}$. Therefore, the formula given for $\mathbf{U}$ in \eqref{l21snfsol} minimizes $\mathscr{L}(\mathbf{U},\mathbf{V}(t))$ at step $t$, that is,
$$\mathscr{L}(\mathbf{U}(t+1),\mathbf{V}(t)) \leq \mathscr{L}(\mathbf{U},\mathbf{V}(t)) \qquad \forall \quad U\in \mathbf{R}^{m \times k}.$$
\end{proof}
Next, let us show when $t\rightarrow \infty$, $\mathbf{V}(t)$ also serves as a good candidate in minimizing $\mathscr{L}(\mathbf{U},\mathbf{V})$ with respect to $\mathbf{V}$ under the non-negativity constraint $\mathbf{V}\geq 0$ via the Karush–Kuhn–Tucker (KKT) condition \cite{Nocedal2006}.
\begin{lemma}
\label{lem:KKT}
Let $\mathbf{V}$ solves the constrained optimization $\displaystyle{\min_{\mathbf{V}\geq 0}}\mathscr{L}(\mathbf{U},\mathbf{V})$. Then $\mathbf{V}$ satisfies the KKT condition
\begin{equation}
\label{kktcondition}
(-\mathbf{DX}^T\mathbf{U} + \mathbf{DVU}^T\mathbf{U} + \alpha \mathbf{LV})_{ij}\mathbf{V}^2_{ij} = 0, \qquad \forall \quad 1\leq i \leq n, \quad 1\leq j \leq k.
\end{equation}
\end{lemma}
\begin{proof}
We define the Lagrangian function $\mathscr{L}_0(\mathbf{V}) = \mathscr{L}(\mathbf{U}, \mathbf{V}) - tr(\mathbf{\Phi} \mathbf{V}^T)$ where the Lagrangian multipliers $\mathbf{\Phi} \in \mathbf{R}_+^{n\times k}$ enforce the non-negativity constraint $\mathbf{V}\geq 0$. The simplified expression \eqref{fhformula} of $\mathscr{L}(\mathbf{U}, \mathbf{V})$ implies
\begin{equation*}
\mathscr{L}_0(\mathbf{V}) = \!tr[\mathbf{XDX}^T]\!-\!2tr[\mathbf{V}^T\mathbf{DX}^T\mathbf{U}]\!+\!tr[\mathbf{U}^T\mathbf{UV}^T\mathbf{DV}]\!+\!\alpha tr[\mathbf{V}^T\mathbf{LV}]\!+\!\beta tr[\mathbf{U\hat{D}}\mathbf{U}^T] \!-\!tr(\mathbf{\Phi} \mathbf{V}^T).
\end{equation*}
Based on KKT Theorem \cite{Nocedal2006}, there hold the zero gradient condition
\begin{equation}
\label{zerogradient}
\frac{\partial \mathscr{L}_0(\mathbf{V})}{\partial \mathbf{V}} = -2\mathbf{DX}^T\mathbf{U} + 2\mathbf{DVU}^T\mathbf{U} + 2\alpha \mathbf{LV} - \mathbf{\Phi} = 0,
\end{equation}
and the complementary slackness condition
\begin{equation}
\label{slackness}
\mathbf{\Phi}_{ij} \mathbf{V}_{ij} = 0, \qquad \forall \quad 1\leq i \leq n, \quad 1\leq j \leq k.
\end{equation}
Combination of \eqref{zerogradient} and \eqref{slackness} yields immediately
\begin{equation*}
(-2\mathbf{DX}^T\mathbf{U} + 2 \mathbf{DVU}^T\mathbf{U} + 2\alpha \mathbf{LV})_{ij}\mathbf{V}_{ij} = 0,
\end{equation*}
which can be rewritten as below, due to the fact that ${V}_{ij} = 0$ is equivalent to ${V}^2_{ij} = 0$:
\begin{equation*}
(-\mathbf{DX}^T\mathbf{U} + \mathbf{DVU}^T\mathbf{U} + \alpha \mathbf{LV})_{ij}\mathbf{V}^2_{ij} = 0, \qquad \forall \quad 1\leq i \leq n, \quad 1\leq j \leq k.
\end{equation*}
\end{proof}
\begin{theorem}
\label{thm:KKT}
The limiting solution of $\mathbf{V}(t)$ when $t\rightarrow \infty$ satisfies the KKT condition \eqref{kktcondition}. Therefore, $\mathbf{V}$ update of \eqref{l21snfsol} is optimal in minimizing $\mathscr{L}(\mathbf{U}, \mathbf{V})$.
\end{theorem}
\begin{proof}
Suppose $\it{V}=\displaystyle{\lim_{t\rightarrow \infty}}\mathbf{V}(t)$. By letting $t\rightarrow \infty$ in the $\mathbf{V}$-update formula \eqref{Vformula} of Lemma~\ref{lem:Vformula}, we have
\begin{equation*}
\it{V}_{ij}\!=\!\it{V}_{ij}\sqrt{\frac{(\mathbf{D}\mathbf{\Phi}^+)_{ij}+(\mathbf{D}\it{V}\mathbf{\Omega}^-)_{ij}+\alpha(\mathbf{W}\it{V})_{ij}}
                {(\mathbf{D}\mathbf{\Phi}^-)_{ij}+(\mathbf{D}\it{V}\mathbf{\Omega}^+)_{ij}+\alpha(\mathbf{\overline{D}}\it{V})_{ij}}} \qquad \forall \quad 1\leq i \leq n, \quad 1\leq j \leq k.
\end{equation*}
By taking square on both sides, then multiplying the equation by the denominator of the right side, we get
\begin{equation*}
\left[(\mathbf{D}\mathbf{\Phi}^-)_{ij}+(\mathbf{D}\it{V}\mathbf{\Omega}^+)_{ij}+\alpha(\mathbf{\overline{D}}\it{V})_{ij}\right] \it{V}^2_{ij}
      = \left[(\mathbf{D}\mathbf{\Phi}^+)_{ij}+(\mathbf{D}\it{V}\mathbf{\Omega}^-)_{ij}+\alpha(\mathbf{W}\it{V})_{ij}\right]\it{V}^2_{ij}.
\end{equation*}
Since $\mathbf{\Omega}^{+}-\mathbf{\Omega}^{-} =\mathbf{\Omega} = \mathbf{U}^T\mathbf{U}$, $\mathbf{\Phi}^{+}-\mathbf{\Phi}^{-}=\mathbf{\Phi}=\mathbf{X}^T\mathbf{U}$ and $\mathbf{L}=\mathbf{\overline{D}}-\mathbf{W}$, the above equation becomes
\begin{equation*}
(-\mathbf{DX}^T\mathbf{U} + \mathbf{D}\it{V}\mathbf{U}^T\mathbf{U} + \alpha \mathbf{L}\it{V})_{ij}\it{V}^2_{ij} = 0.
\end{equation*}
Therefore, the limiting solution $\it{V}$ satisfies the KKT condition \eqref{kktcondition}.
\end{proof}
Finally, let us show that $\mathscr{L}(\mathbf{U}(t),\mathbf{V}(t)$ and $\mathscr{J}(\mathbf{U}(t),\mathbf{V}(t))$ go towards each other when $t\rightarrow \infty$.
\begin{theorem}
\label{thm:LJrelation}
Let $(\it{U},\it{V})=\displaystyle{\lim_{t\rightarrow \infty}}(\mathbf{U}(t),\mathbf{V}(t))$, then $\mathscr{L}(\it{U},\it{V})=\mathscr{J}(\it{U},\it{V})$.
\end{theorem}
\begin{proof}
Based on the \textit{proxy loss} function definition in \eqref{proxyloss}, the limit solution $(\it{U}, \it{V})$ satisfies
\begin{equation}
\label{limitloss}
\mathscr{L}(\it{U},\it{V})\!=\!tr[(\mathbf{X}\!-\!\it{UV}^T)\mathbf{D}(\mathbf{X}\!-\!\it{UV}^T)^{T}]
+\alpha tr[\it{V}^T\mathbf{L}\it{V}]
+\beta tr[\it{U}\mathbf{\hat{D}}\it{U}^T]
\end{equation}
where $\mathbf{D}_{ii} \!=\! \frac{1}{\|\mathbf{X}_i-\it{U}\it{V}_i^T\|_{2}}$, $\mathbf{\hat{D}_{ii}} \!=\! \frac{1}{\|\it{U}_i\|_{2}}$,
$\mathbf{L} \!=\! \mathbf{\overline{D}}\!-\!\mathbf{W}$, $\mathbf{W}_{ij}=\frac{w_{ij}}{\|\it{V}_i-\it{V}_j\|_{2}} \enspace i\neq j$,
$\mathbf{\overline{D}}_{ii}\!=\!\displaystyle{\sum_{j\neq i}}\mathbf{W}_{ij}$. Meanwhile,
\begin{equation}
\label{limitobj}
\mathscr{J}(\it{U},\it{V}) = \|\mathbf{X}-\it{UV}^T\|_{2,1} +\alpha \sum_{i<j} \|\it{V}_i-\it{V}_j\|_{2}w_{ij} + \beta\|\it{U}\|_{2,1}.
\end{equation}
Let $C_i=\mathbf{X}_i\!-\!\it{UV}_i^T$ denote the $i$-th column of $\mathbf{X}\!-\!\it{UV}^T$, then the first term of $\mathscr{L}(\it{U},\it{V})$ in \eqref{limitloss} becomes
\begin{equation}
tr[(\mathbf{X}\!-\!\mathbf{\it UV}^T)\mathbf{D}(\mathbf{X}\!-\!\mathbf{\it UV}^T)^{T}] \!=\! tr\left[\sum_{i=1}^n \frac{1}{\|C_i\|_2}C_iC_i^T\right]
   \!=\! \sum_{i=1}^n \frac{1}{\|C_i\|_2} \|C_i\|^2_2 \!=\! \sum_{i=1}^n \|C_i\|_2 \!=\! \|\mathbf{X}\!-\!\it{UV}^T\|_{2,1}.
\label{limitloss1}
\end{equation}
Similarly, by assuming $U_i$ as the $i$-th column of $\it{U}$, the third term of $\mathscr{L}(\it{U},\it{V})$ in \eqref{limitloss} becomes
\begin{equation}
\beta tr[\it{U}\mathbf{\hat{D}}\it{U}^T] \!=\! \beta tr\left[\sum_{i=1}^k \frac{1}{\|U_i\|_2}U_iU_i^T\right]
   \!=\! \beta \sum_{i=1}^k \frac{1}{\|U_i\|_2} \|U_i\|^2_2 \!=\! \beta \sum_{i=1}^k \|U_i\|_2 \!=\! \beta \|\it{U}\|_{2,1}.
\label{limitloss3}
\end{equation}
Finally, by assuming $V_i$ as the $i$-th row of $\it{V}$ and $\langle\it{V}_i,\it{V}_j\rangle$ as the vector inner product, the second term of $\mathscr{L}(\it{U},\it{V})$ in \eqref{limitloss} becomes
\begin{align}
\alpha tr[\it{V}^T\mathbf{L}\it{V}] & = \alpha tr[\it{V}^T\mathbf{\overline{D}}\it{V}] - \alpha tr[\it{V}^T\mathbf{W}\it{V}]
        = \alpha \sum_{i=1}^n \mathbf{\overline{D}}_{ii} \langle\it{V}_i,\it{V}_i\rangle
          -\alpha \sum_{i\neq j} \mathbf{W}_{ij} \langle\it{V}_i,\it{V}_j\rangle  \nonumber \\
      & = \alpha \sum_{i=1}^n \left(\sum_{j\neq i} \mathbf{W}_{ij}\right) \langle\it{V}_i,\it{V}_i\rangle
          -\alpha \sum_{i\neq j} \mathbf{W}_{ij} \langle\it{V}_i,\it{V}_j\rangle  \nonumber \\
      & = \alpha \sum_{i<j} \mathbf{W}_{ij} \left[\langle\it{V}_i,\it{V}_i\rangle-2\langle\it{V}_i,\it{V}_j\rangle
                                                                    +\langle\it{V}_j,\it{V}_j\rangle \right] \nonumber \\
      & = \alpha \sum_{i<j} \frac{w_{ij}}{\|\it{V}_i-\it{V}_j\|_{2}} \|\it{V}_i-\it{V}_j\|_{2}^2 \nonumber \\
      & = \alpha \sum_{i<j} \|\it{V}_i-\it{V}_j\|_{2}w_{ij}.
\label{limitloss2}
\end{align}
Summing up both sides of \eqref{limitloss1}, \eqref{limitloss3} and \eqref{limitloss2} yields \eqref{limitloss} and \eqref{limitobj}, respectively. Therefore,
$$\mathscr{L}(\it{U},\it{V})=\mathscr{J}(\it{U},\it{V}).$$
\end{proof}
In summary, for the sequence $\{\mathbf{U(t), V(t)}\}$ constructed by \eqref{l21snfsol}, when $t\rightarrow \infty$, the \textit{proxy loss} function $\mathscr{L}(\mathbf{U},\mathbf{V})$ reaches its minima first towards $\mathbf{U}$ with fixed $\mathbf{V}$, then towards $\mathbf{V}$ with fixed $\mathbf{U}$ at each step. Meanwhile, $\mathscr{L}(\mathbf{U},\mathbf{V})$ and the objective function $\mathscr{J}(\mathbf{U},\mathbf{V})$ go towards each other when $t\rightarrow \infty$. Therefore, $\{\mathbf{U(t), V(t)}\}$ serves as a solid candidate in minimizing $\mathscr{J}(\mathbf{U},\mathbf{V})$ iteratively. The following numerical experiments verify robustness of $L_{2,1}$ SNF algorithm \eqref{l21snfsol}.

\section{Experimental Results}
\label{sec:experiment}

In this section, we compare the performance of the original SNF algorithm, graph regularized SNF algorithm, and $L_{2,1}$ SNF algorithm by numerical experiments on three benchmark mixed-sign datasets  and several randomized mixed-sign matrices.

\subsection{Experimental Setting}
\label{subsec:setting}

In our experiments, we use three benchmark mixed-sign datasets which are generated from real world data collection and are widely used in the clustering literature. Table~\ref{dataset} summarizes their statistics and provides their URLs for detailed information.

Ionosphere: It is from the UCI repository and consists of radar data collected in Goose Bay, Canada. It has $2$ classes with $351$ instances and each instance has $34$ numeric attributes ($17$ pulse complex numbers).

Waveform: It is also from the UCI repository and is a collection of simulated time series signals for signal processing applications. It has $3$ classes with $5000$ instances and each instance has $21$ numerical attributes.

USPST: Ir is the test split of the USPS system, and each image is presented at the resolution of $16\times 16$ pixels. The dataset has $10$ classes representing digits $0,1,2,\cdots,9$ with $2007$ images. It truly represents the original USPS system with 7291 images and meanwhile maintains similar computational load as the other two datasets.

\begin{table}[htb]
\centering
\renewcommand{\arraystretch}{0.90}
\begin{tabular}{|c|c|c|c|c|}
\hline
Dataset & Instance & Attribute & Classes & URL \\
\hline
Ionosphere &  351  &  34   &  2  & \href{https://archive.ics.uci.edu/dataset/52/ionosphere}{UCI} \\
Waveform   &  5000 &  21   &  3  & \href{https://archive.ics.uci.edu/dataset/107/waveform+database+generator+version+1}{UCI} \\
USPST      &  2007 &  256  &  10 & \href{https://git-disl.github.io/GTDLBench/datasets/usps_dataset}{Github} \\
\hline
\end{tabular}
\caption{Statistics of the Benchmark Datasets}
\label{dataset}
\end{table}

To simplify our experiments, we set the number $p$ of nearest neighbors for data graph to an near-optimal value $5$ for GR SNF and $L_{2,1}$ SNF algorithms, each of which has $2$ remaining tuning parameters $\alpha$ and $\beta$. To compare these methods fairly, we perform grid search in the parameter spaces for $\alpha$ and $\beta$ as a pair and report the best results for them. There is no parameter selection for SNF algorithm.

The clustering result is evaluated by comparing the obtained label of each sample with the true label provided by the dataset. To capture the obtained labels for all the vectors in $\mathbf{X}$ after decomposition $\mathbf{X}\approx \mathbf{UV}^T$, we apply k-means clustering method with $k$ clusters towards all the rows of $\mathbf{V}$, which represents the coordinates of basis vectors in $\mathbf{U}$ for the approximation of vectors in $\mathbf{X}$. Then the vectors within each cluster are assigned a common label that appears the most among the true labels of these vectors. We adopt two widely-used metrics, accuracy of clustering and normalized mutual information towards the true and obtained labels of the dataset to measure the clustering performance. These metrics are more reasonable for algorithm comparison than the traditional error estimate of $\mathbf{X}-\mathbf{UV}^T$ as different algorithms may choose different values of $\alpha$ and $\beta$ to minimize their own objective functions.

{\bf Accuracy of clustering (ACC)} is defined as $\textbf{ACC}(\mathbf{Y},\mathbf{C}) = \frac{\sum_{i=1}^n \delta(y_i,c_i)}{n}$,
where $\mathbf{Y}=(y_1, \cdots, y_n)$ are the true class labels, and $\mathbf{C}=(c_1, \cdots, c_n)$ are the obtained cluster labels of the given dataset $\mathbf{X}$. $\delta(x,y)$ is the delta function where $\delta(x,y)=1$ if $x=y$ and $\delta(x,y)=0$ otherwise. A large ACC value indicates a better clustering performance.

{\bf Normalized mutual information (NMI)} is defined as $\textbf{NMI}(\mathbf{Y},\mathbf{C})=\frac{I(\mathbf{Y},\mathbf{C})}
{\max(H(\mathbf{Y}),H(\mathbf{C}))}$, where $I(\mathbf{Y},\mathbf{C})=\sum_{y\in Y, c\in C}p(y,c)\log_2\frac{p(y,c)}{p(y)p(c)}$ is the mutual information between true class labels
$\mathbf{Y}$ and clustering labels $\mathbf{C}$, $H(\mathbf{Y})=-\sum_{y\in Y}p(y)\log_2p(y)$ and $H(\mathbf{C})=-\sum_{c\in C}p(c)\log_2p(c)$ are the entropies for $\mathbf{Y}$ and
$\mathbf{C}$, where $p(y)$ is the probability that $y\in \mathbf{Y}$ is selected, $p(c)$ is the probability that $c\in \mathbf{C}$ is selected, and $p(y,c)$ is the joint probability that $y\in \mathbf{Y}$ and $c\in \mathbf{C}$ are selected simultaneously. A larger NMI value also indicates a better clustering performance.

Since the class numbers of Ionosphere, Waveform and USPST are $2$, $3$ and $10$, respectively, we set the cluster number $k$ larger than the class number for each dataset. To be more specific, we set $k = 4, 5, 6, 7$ for Ionosphere; $k = 8, 10, 12, 14$ for Waveform; $k = 12, 16, 20, 24$ for USPST. For each given cluster number, $20$ test runs are conducted on different randomly chosen data which account for $90\%$ of the whole dataset. The number of iterations for each test is $500$. The initialization values of $\mathbf{U}$ and $\mathbf{V}$ are chosen randomly in the range of $[-1, 1]$ and $[0, 1]$, respectively. The average of the mean and standard deviation for ACC and NMI from those tests are reported for each algorithm.

\subsection{Performance Comparison of All Algorithms}
\label{subsec:performance}

\begin{table}[htb]
\hspace{5mm}
\renewcommand{\arraystretch}{0.90}
\begin{tabular}{|c|c|c|c||c|c|c|c|}
\hline
\multicolumn{4}{|c||}{ACC (\%) Mean$\pm$SD} & \multicolumn{3}{|c|}{NMI (\%) Mean$\pm$SD}\\
\hline
$k$ & SNF & GR SNF & $L_{2,1}$ SNF & SNF & GR SNF & $L_{2,1}$ SNF \\
\hline
4 & 82.40$\pm$2.35 & 85.16$\pm$1.21 & 85.24$\pm$1.49 & 33.28$\pm$6.31 & 36.83$\pm$3.24 & 37.24$\pm$3.88 \\
5 & 82.04$\pm$3.13 & 85.08$\pm$0.90 & 85.65$\pm$0.57 & 32.21$\pm$7.27 & 36.69$\pm$2.51 & 38.43$\pm$1.72 \\
6 & 81.59$\pm$1.85 & 85.49$\pm$0.89 & 85.60$\pm$0.71 & 29.91$\pm$4.06 & 38.09$\pm$2.48 & 38.34$\pm$2.08 \\
7 & 81.98$\pm$2.36 & 84.81$\pm$0.66 & 85.33$\pm$0.60 & 31.56$\pm$5.11 & 36.55$\pm$1.93 & 37.44$\pm$1.56 \\
\hline
\end{tabular}
\caption{Clustering Performance on Ionosphere}
\label{perf_ionosphere}
\end{table}

\begin{table}[htb]
\hspace{5mm}
\renewcommand{\arraystretch}{0.9}
\begin{tabular}{|c|c|c|c||c|c|c|c|}
\hline
\multicolumn{4}{|c||}{ACC (\%) Mean$\pm$SD} & \multicolumn{3}{|c|}{NMI (\%) Mean$\pm$SD}\\
\hline
$k$ & SNF & GR SNF & $L_{2,1}$ SNF & SNF & GR SNF & $L_{2,1}$ SNF \\
\hline
8  & 51.03$\pm$2.34 & 74.86$\pm$2.18 & 77.98$\pm$0.57 & 11.15$\pm$4.20 & 40.48$\pm$4.37 & 47.13$\pm$0.98 \\
10 & 50.94$\pm$2.84 & 74.21$\pm$0.57 & 81.22$\pm$0.36 & 9.88$\pm$4.18 & 39.30$\pm$2.94 & 50.26$\pm$0.75 \\
12 & 47.83$\pm$2.77 & 75.17$\pm$0.96 & 81.45$\pm$0.47 & 6.80$\pm$3.18 & 40.20$\pm$3.26 & 49.79$\pm$1.30 \\
14 & 48.29$\pm$1.91 & 77.50$\pm$0.96 & 80.65$\pm$0.79 & 6.83$\pm$2.27 & 42.59$\pm$2.25 & 46.86$\pm$2.03 \\
\hline
\end{tabular}
\caption{Clustering Performance on Waveform}
\label{perf_waveform}
\end{table}

\begin{table}[htb]
\hspace{5mm}
\renewcommand{\arraystretch}{0.9}
\begin{tabular}{|c|c|c|c||c|c|c|c|}
\hline
\multicolumn{4}{|c||}{ACC (\%) Mean$\pm$SD} & \multicolumn{3}{|c|}{NMI (\%) Mean$\pm$SD}\\
\hline
$k$ & SNF & GR SNF & $L_{2,1}$ SNF & SNF & GR SNF & $L_{2,1}$ SNF \\
\hline
12 & 67.77$\pm$3.75 & 77.74$\pm$1.22 & 80.49$\pm$1.81 & 55.61$\pm$2.92 & 66.02$\pm$0.86 & 71.10$\pm$1.20 \\
16 & 68.41$\pm$2.26 & 77.62$\pm$0.79 & 81.55$\pm$1.08 & 55.17$\pm$2.24 & 67.76$\pm$0.97 & 72.33$\pm$0.92 \\
20 & 69.39$\pm$2.41 & 80.02$\pm$0.99 & 82.56$\pm$1.15 & 55.55$\pm$2.55 & 69.75$\pm$1.00 & 73.04$\pm$0.95 \\
24 & 71.74$\pm$1.85 & 81.96$\pm$0.89 & 83.94$\pm$1.15 & 56.74$\pm$2.11 & 70.97$\pm$1.05 & 73.70$\pm$1.35 \\
\hline
\end{tabular}
\caption{Clustering Performance on USPST}
\label{perf_uspst}
\end{table}

Table~\ref{perf_ionosphere}, Table~\ref{perf_waveform} and Table~\ref{perf_uspst} show the clustering results on Ionosphere, Waveform and USPST datasets, respectively. The mean and standard deviation (SD) of ACC and NMI are reported in the tables. The optimal parameter values chosen for $L_{2,1}$ SNF are $\alpha=0.1, \beta=2.25$ for Ionosphere, $\alpha=0.1, \beta=100$ for Waveform, and $\alpha=1, \beta=15$ for USPST, respectively.

Based on these experiments, we can see that $L_{2,1}$ SNF outperforms SNF and GR SNF algorithms across each benchmark dataset. Next, we demonstrate similar performance of these algorithms on the same datasets with Gaussian noise at different levels.

\subsection{Noise Handling Performance of All Algorithms}
\label{subsec:noise}

In this subsection, we compare the clustering performance of all the algorithms on the above three datasets under the influence of Gaussian noise. The Gaussian noise usually happens in amplifiers or detectors and generates disturbs in the gray values of the given data. It can be described as below:
\begin{equation}
 \begin{aligned}
    P(g) = \sqrt{\frac{1}{2\pi \sigma^2}}e^{-\frac{(g - \mu)^2}{2 \sigma ^2}}
 \end{aligned}
 \label{gaussian_noise}
\end{equation}
where $g$ indicates the gray value, $\mu$ is the mean, and $\sigma$ is the standard deviation. The Gaussian function has its peak at the mean, and its “spread” increases with the standard deviation so that the function reaches $0.607$ times its maximum at $g+\sigma$ and $g-\sigma$. Since all three datasets have mixed-sign data around zero, we simply set the mean $\mu=0$. We then choose different values of $\sigma$ to reflect different noise levels for our experiments.

\begin{figure}[htb]
\hspace{-5mm}
\subfloat{\includegraphics[width=0.55\textwidth]{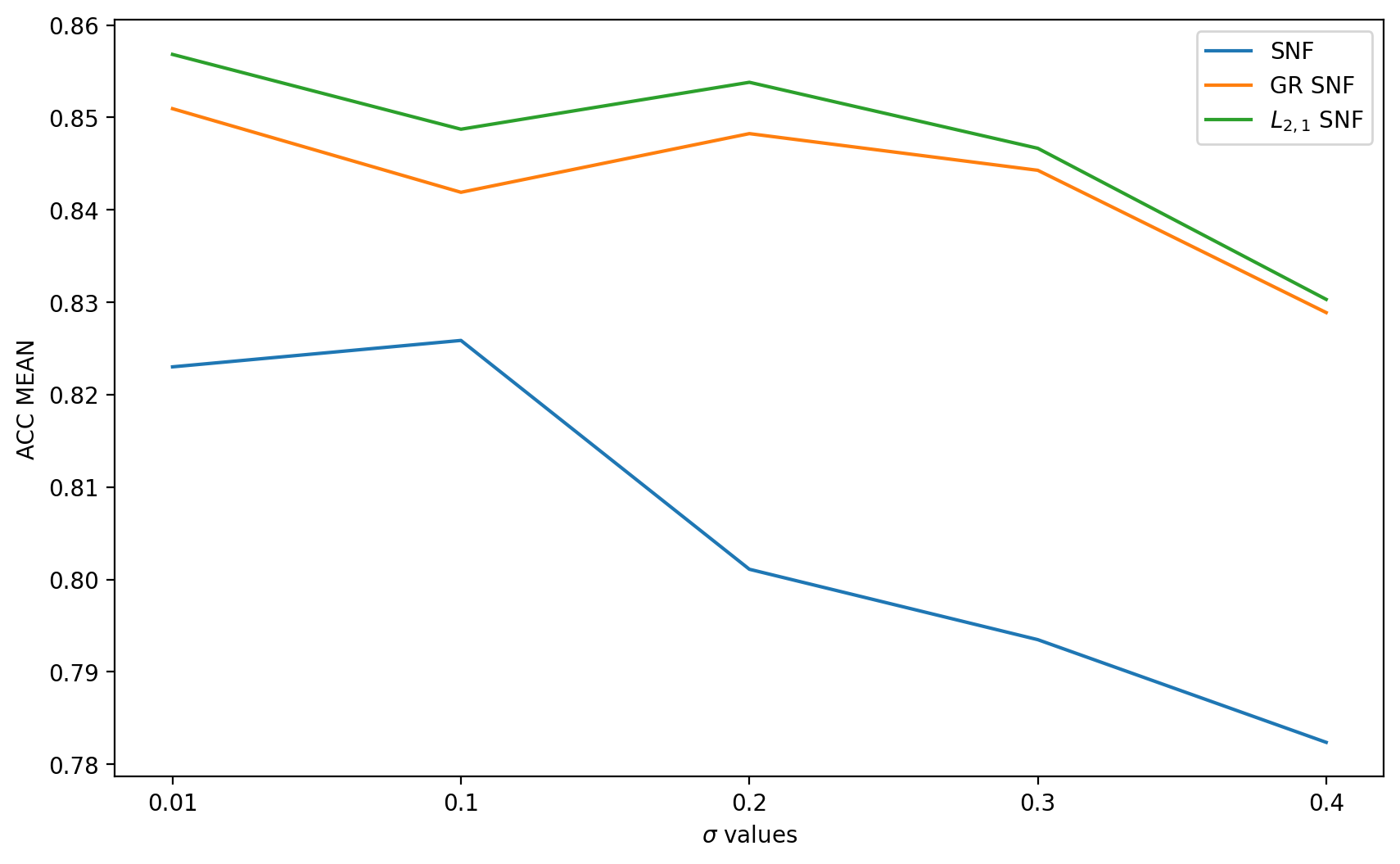}}
\subfloat{\includegraphics[width=0.55\textwidth]{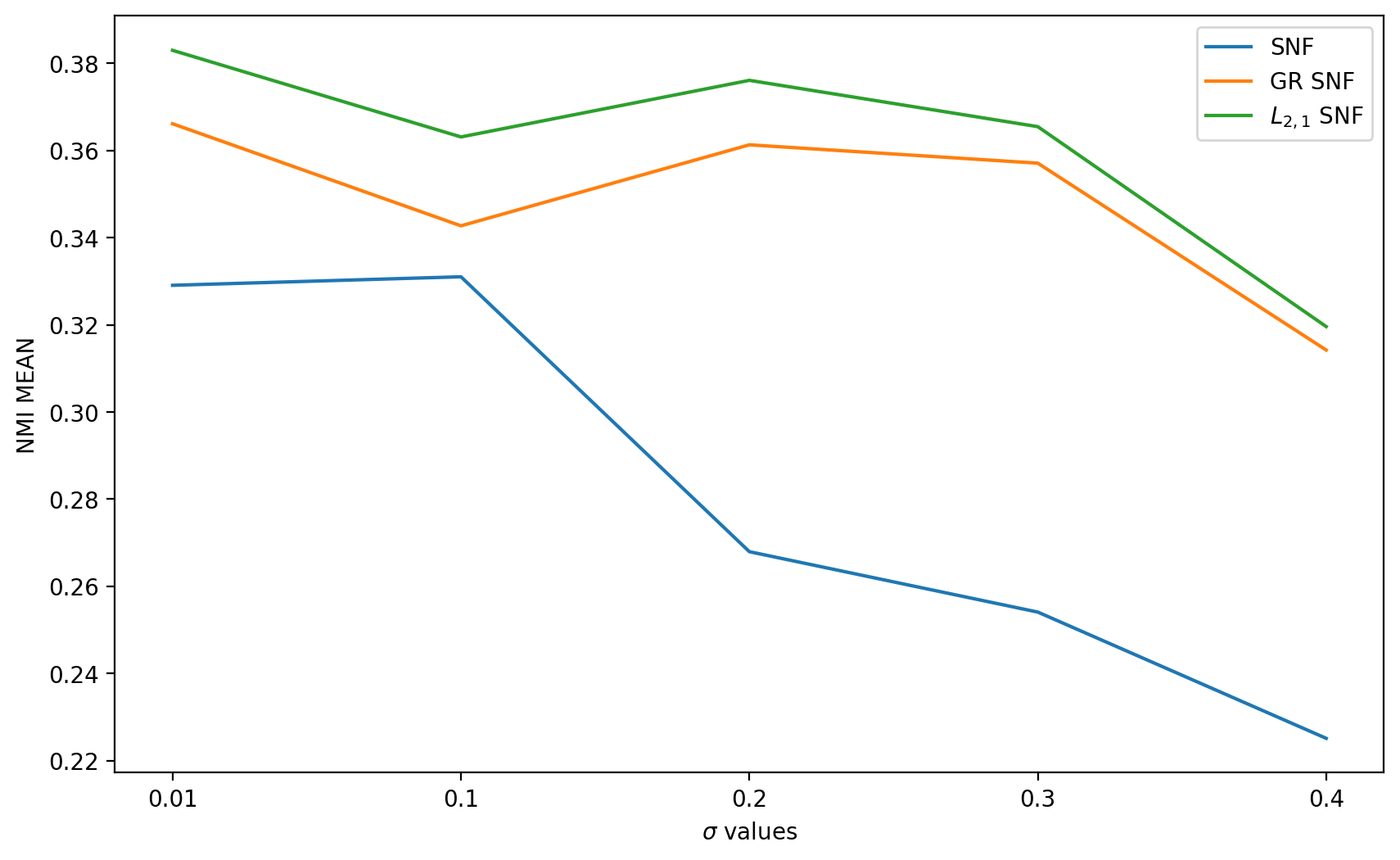}}
\caption{Clustering Performance with $\sigma$ noise level on Ionosphere.}
\label{noise_ionosphere}
\end{figure}

\begin{figure}[htb]
\hspace{-5mm}
\subfloat{\includegraphics[width=0.55\textwidth]{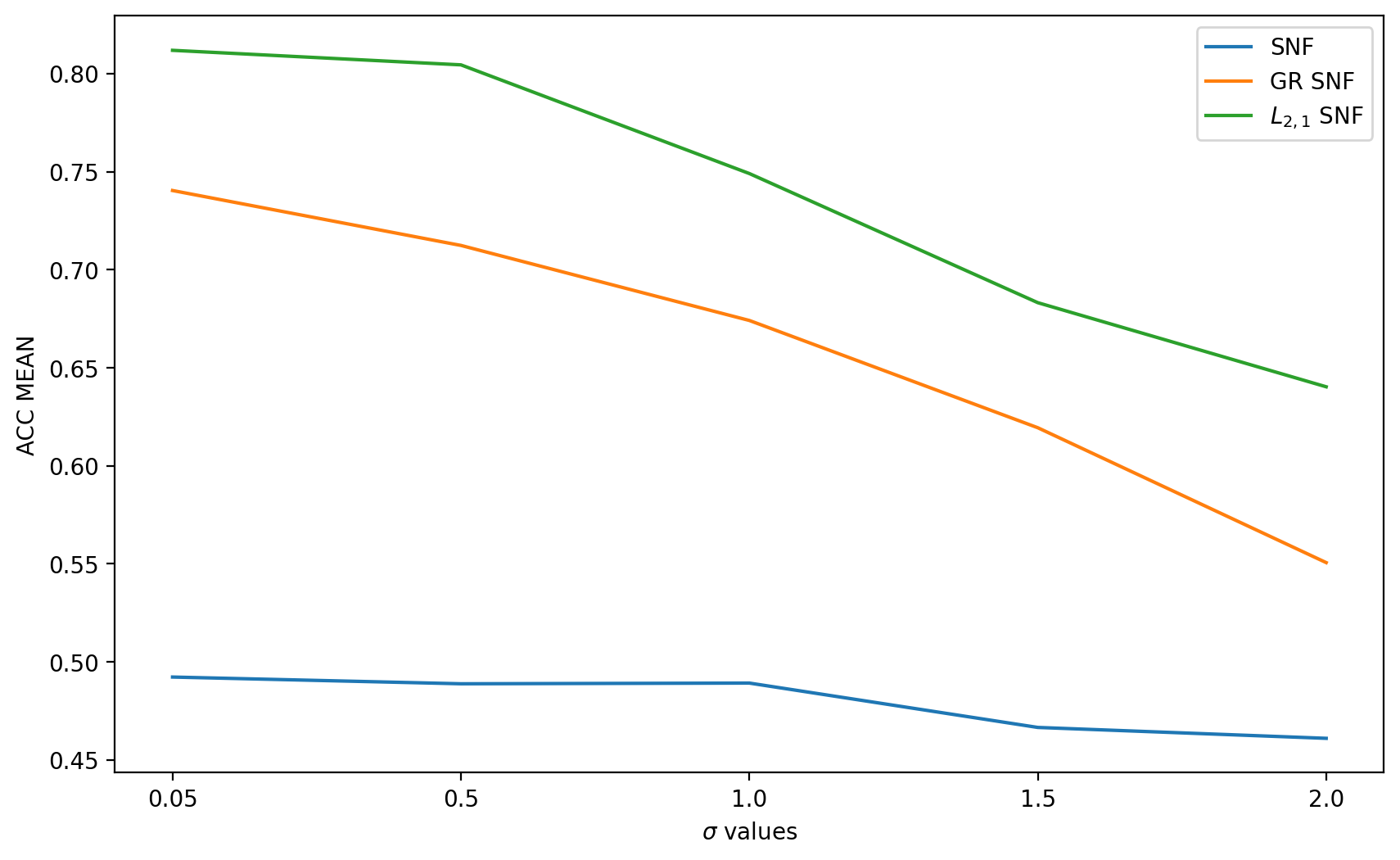}}
\subfloat{\includegraphics[width=0.55\textwidth]{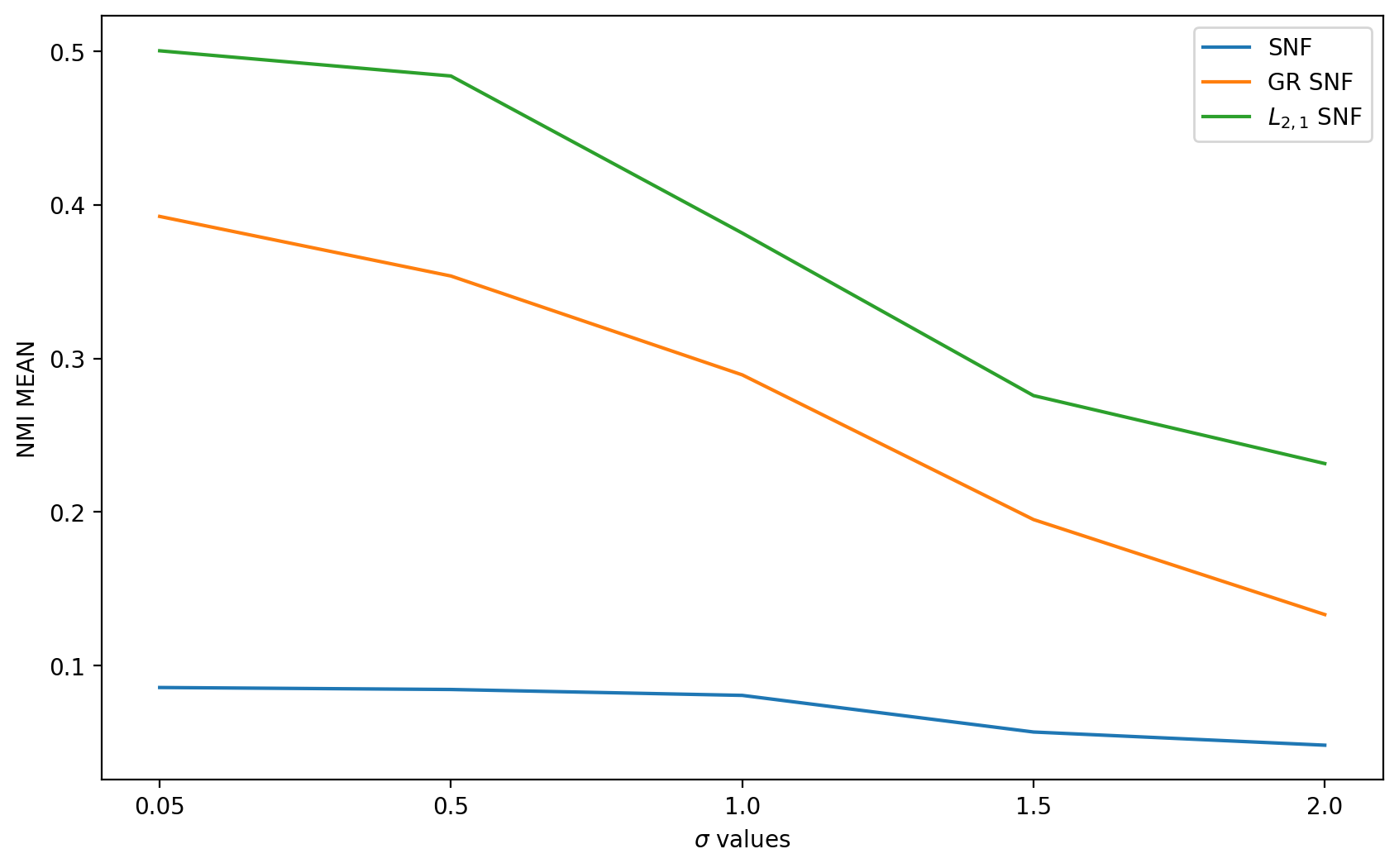}}
\caption{Clustering Performance with $\sigma$ noise level on Waveform.}
\label{noise_waveform}
\end{figure}

\begin{figure}[htb]
\hspace{-5mm}
\subfloat{\includegraphics[width=0.55\textwidth]{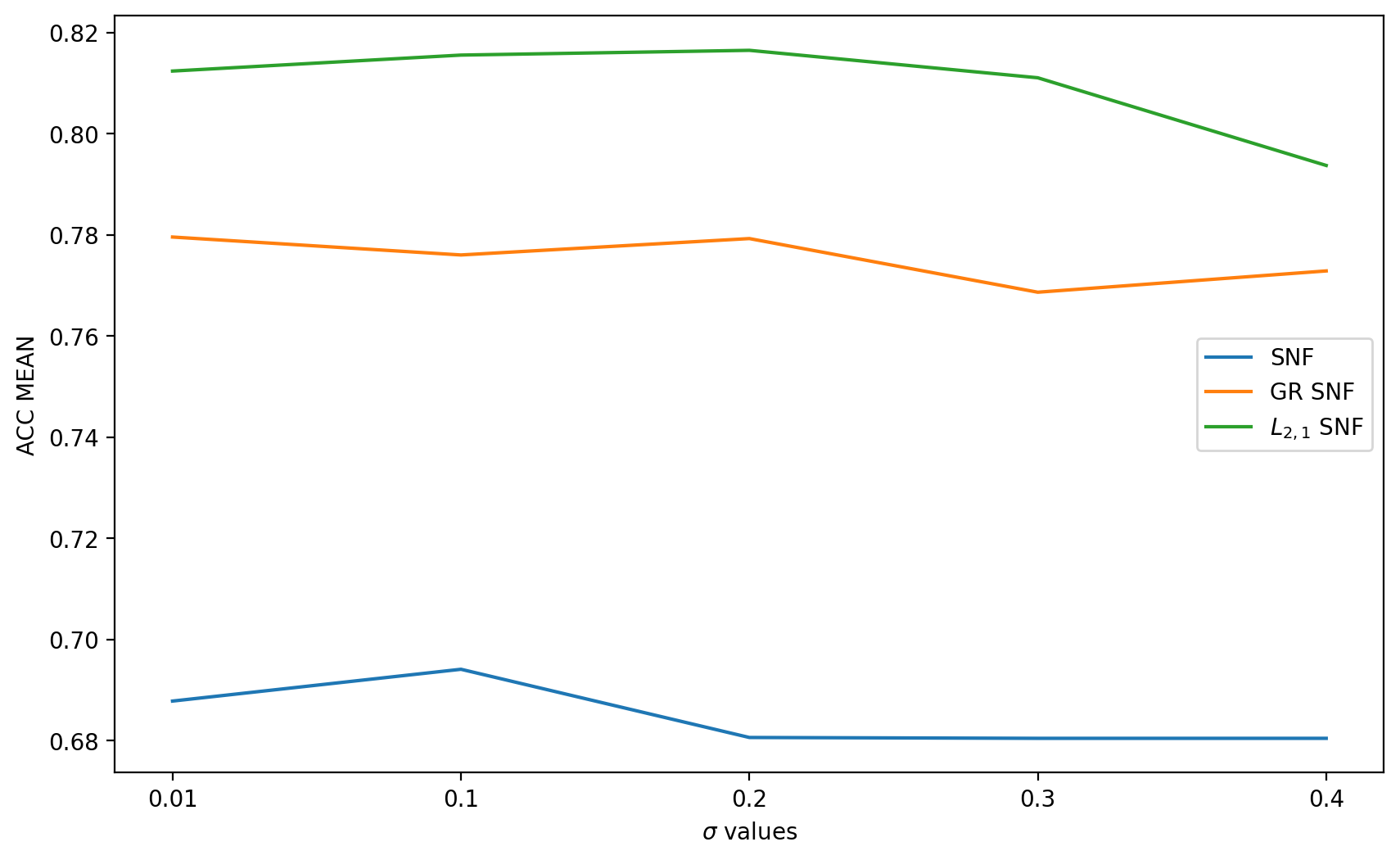}}
\subfloat{\includegraphics[width=0.55\textwidth]{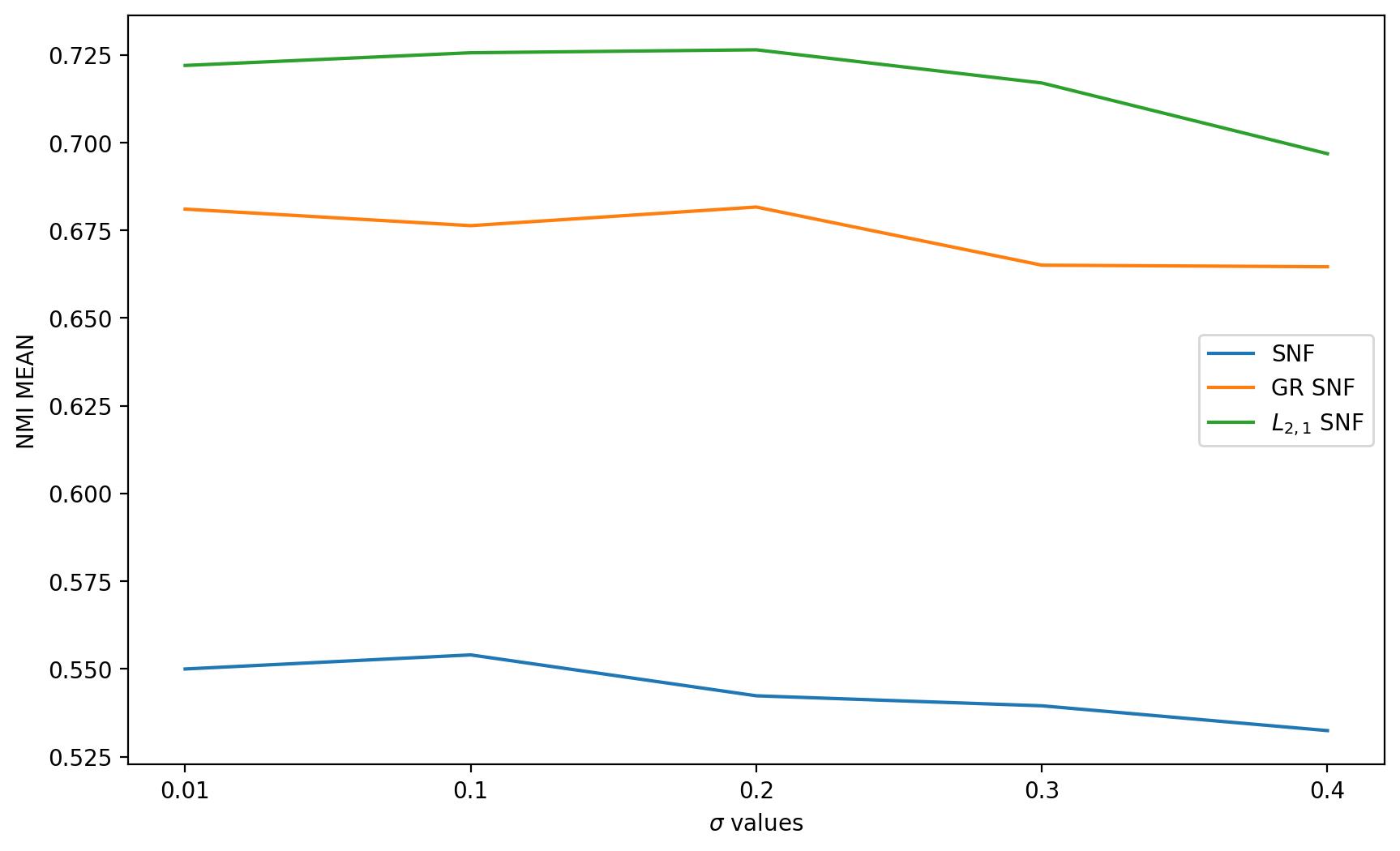}}
\caption{Clustering Performance with $\sigma$ noise level on USPST.}
\label{noise_uspst}
\end{figure}

Figure~\ref{noise_ionosphere}, Figure~\ref{noise_waveform} and Figure~\ref{noise_uspst} show the clustering results on Ionosphere, Waveform and USPST datasets, respectively. The ACC and NMI means are reported in the figures where the horizontal axis denotes the values of $\sigma$ as noise levels. The corresponding ACC and NMI standard deviation values are very similar to those without noise as reported in Tables~\ref{perf_ionosphere},~\ref{perf_waveform},~\ref{perf_uspst} and are thus ignored for brevity.

Based on these experiments, it is observed that for all the datasets, the performance of all SNF algorithms deteriorate as the noise level increases. However, $L_{2,1}$ SNF still outperforms SNF and GR SNF across all the datasets at different levels of Gaussian noise.

\begin{figure}[htb]
\hspace{-5mm}
\includegraphics[width=1.1\textwidth]{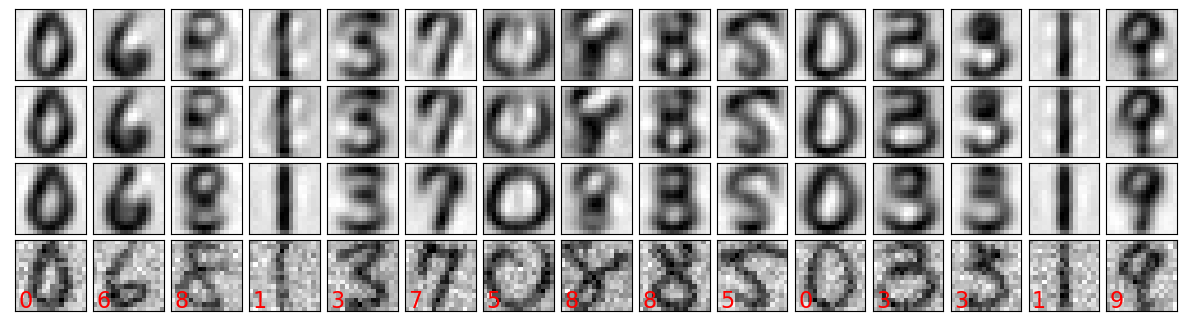}
\caption{Images of $15$ randomly selected numbers reconstructed by all SNF algorithms on USPST. Top: SNF; Second from Top: GR SNF; Second from Bottom: $L_{2,1}$ SNF; Bottom: Original image with noise where the true label is shown in red at the lower left corner.}
\label{image_uspst}
\end{figure}

As USPST dataset is an image collection of numbers $0, 1, 2, \cdots, 9$, we run an experiment with cluster number $k=16$ on USPST dataset which carries Gaussian noise with $\sigma=0.30$. The ACC and NMI means are reported as $(0.657, 0.528)$, $(0.791, 0.68)$ and $(0.822, 0.735)$ for SNF, GR SNF and $L_{2,1}$ SNF, respectively. Meanwhile, Figure~\ref{image_uspst} shows the images of $15$ randomly selected numbers reconstructed by all the algorithms as well as the original image with given Gaussian noise. It is observed that the images by $L_{2,1}$ SNF achieve higher precision to recapture the original images than those by the other two SNF algorithms.

\subsection{Parameter Investigation for $L_{2,1}$ SNF}
\label{subsec:parameter}

$L_{2,1}$ SNF has two essential parameters: $\alpha$ measures the weight of the graph Laplacian; $\beta$ controls the sparse degree of the basis matrix. We investigate their influence on the ACC and NMI means for each dataset by varying $\alpha$ and $\beta$ as a pair. For fairness of comparison, we choose $\alpha$ from a wide feasible range $\{10^{-3}, 10^{-2}, 10^{-1}, 1, 10, 10^2\}$ commonly adopted in the literature and $\beta$ from $\{10^{-3}, 10^{-2}, 10^{-1}, 1, 10, 10^2, 10^3\}$. Once an $(\alpha, \beta)$ pair is identified to achieve the best ACC and NMI means, we then conduct a refined search in the neighborhood of that pair until a new pair with better values can be found. We run the experiments with the cluster number set to $k=5$ for Ionosphere, $k=10$ for Waveform, and $k=16$ for USPST, respectively.

Table~\ref{ab_ionosphere}, Table~\ref{ab_waveform}, and Table~\ref{ab_uspst} show how the ACC and NMI mean values vary with $(\alpha,\beta)$ for each dataset. The optimal $(\alpha, \beta)$ pair for each dataset is marked in red in each table, respectively: $(0.1, 2.25)$ for Ionosphere, $(0.1, 100)$ for Waveform, and $(1.0, 15.0)$ for USPST. It is observed that all three datasets follow the same pattern: when the $(\alpha, \beta)$ pair approaches towards the optimal pair from any direction, the performance of $L_{2,1}$ SNF improves gradually; when the pair moves away from the optimal pair, the performance of the algorithm deteriorates. Therefore, when assigned appropriate values, the graph Laplacian and the sparseness constraint are surely helpful for a better data representation.

\begin{table}[htb]
\hspace{5mm}
\renewcommand{\arraystretch}{0.9}
\begin{tabular}{|p{0.9cm}|p{1.8cm}|p{1.8cm}|p{1.8cm}|p{1.8cm}|p{1.8cm}|p{1.8cm}|}
\hline
\diagbox{$\beta$}{$\alpha$} & 0.001 & 0.01 & 0.1 & 1.0 & 10 & 100 \\
\hline
0.001  & (80.3,26.6) & (82.9,31.5) & (85.1,39.0) & (83.4,33.9) & (78.1,22.6) & (76.6,19.8) \\
\hline
0.01   & (80.5,27.8) & (83.3,33.0) & (84.5,37.5) & (83.6,34.0) & (79.5,25.5) & (78.3,22.9) \\
\hline
0.1    & (80.8,27.9) & (83.2,32.5) & (84.7,36.9) & (84.4,36.1) & (79.3,24.6) & (75.2,18.7) \\
\hline
1.0    & (77.8,22.2) & (83.5,33.3) & (85.2,36.8) & (83.5,32.8) & (79.8,24.1) & (75.9,20.2) \\
\hline
2.25   & (75.6,18.7) & (81.2,29.0) & \textcolor{red}{(85.6,38.2)} & (81.4,28.2) & (79.1,23.9) & (77.3,22.8) \\
\hline
10   & (73.2,14.4) & (84.0,34.0) & (82.8,31.0) & (80.2,26.6) & (76.8,20.2) & (75.9,19.8) \\
\hline
100  & (82.3,31.5) & (69.0,9.1) & (69.2,9.7) & (78.8,23.6) & (76.8,21.3) & (75.3,18.5) \\
\hline
1000 & (75.5,18.4) & (79.8,25.6) & (79.0,24.5) & (77.1,21.4) & (76.2,20.2) & (76.0,19.5) \\
\hline
\end{tabular}
\caption{(ACC, NMI) Means by $L_{2,1}$ SNF with $(\alpha, \beta)$ on Ionosphere}
\label{ab_ionosphere}
\end{table}

\begin{table}[htb]
\hspace{5mm}
\renewcommand{\arraystretch}{0.9}
\begin{tabular}{|p{0.9cm}|p{1.8cm}|p{1.8cm}|p{1.8cm}|p{1.8cm}|p{1.8cm}|p{1.8cm}|}
\hline
\diagbox{$\beta$}{$\alpha$} & 0.001 & 0.01 & 0.1 & 1.0 & 10 & 100 \\
\hline
0.001  & (50.8,10.4) & (54.7,13.8) & (72.8,36.7) & (80.2,49.0) & (77.8,46.7) & (71.7,38.6) \\
\hline
0.01   & (52.0,9.8) & (55.0,14.7) & (72.7,37.1) & (80.5,49.3) & (77.5,45.8) & (72.0,39.0) \\
\hline
0.1    & (52.0,12.5) & (54.1,14.0) & (72.9,37.1) & (79.8,47.5) & (77.1,47.3) & (71.5,38.9) \\
\hline
1.0    & (54.1,12.2) & (60.9,18.8) & (73.7,37.3) & (74.6,41.0) & (77.6,47.1) & (71.1,38.2) \\
\hline
10   & (63.8,24.7) & (69.0,31.8) & (79.1,47.8) & (80.9,48.8) & (76.8,46.1) & (70.9,37.5) \\
\hline
100  & (55.6,13.1) & (69.7,31.2) & \textcolor{red}{(81.2,50.3)} & (78.5,47.8) & (71.4,38.4) & (70.0,35.6) \\
\hline
1000 & (72.8,37.1) & (80.9,49.8) & (76.5,43.8) & (62.6,23.7) & (69.3,35.6) & (70.5,35.6) \\
\hline
\end{tabular}
\caption{(ACC, NMI) Means by $L_{2,1}$ SNF with $(\alpha, \beta)$ on Waveform}
\label{ab_waveform}
\end{table}

\begin{table}[htb]
\hspace{5mm}
\renewcommand{\arraystretch}{0.9}
\begin{tabular}{|p{0.9cm}|p{1.8cm}|p{1.8cm}|p{1.8cm}|p{1.8cm}|p{1.8cm}|p{1.8cm}|}
\hline
\diagbox{$\beta$}{$\alpha$} & 0.001 & 0.01 & 0.1 & 1.0 & 10 & 100 \\
\hline
0.001  & (67.7,52.9) & (68.5,55.5) & (74.7,64.3) & (69.5,61.7) & (74.5,69.9) & (60.6,51.3) \\
\hline
0.01   & (66.5,52.2) & (69.5,56.0) & (75.3,64.7) & (67.8,60.1) & (77.1,69.9) & (61.4,51.7) \\
\hline
0.1    & (68.7,55.2) & (70.8,58.2) & (74.6,64.6) & (65.1,59.0) & (77.0,70.0) & (59.8,50.8) \\
\hline
1.0    & (74.3,63.3) & (75.1,63.9) & (77.7,67.2) & (77.7,68.2) & (80.6,71.7) & (61.2,51.8) \\
\hline
10    & (74.5,62.5) & (75.6,64.5) & (78.4,68.2) & (79.9,70.9) & (70.2,62.7) & (57.9,48.6) \\
\hline
15   & (73.4,61.6) & (76.9,65.7) & (78.8,68.6) & \textcolor{red}{(81.7,72.4)} & (63.6,55.9) & (54.8,45.1) \\
\hline
100  & (76.0,64.4) & (78.7,67.8) & (78.7,68.9) & (55.2,46.2) & (56.0,46.3) & (53.7,43.9) \\
\hline
1000 & (76.9,66.0) & (61.8,54.8) & (37.3,21.9) & (59.9,50.9) & (54.2,44.1) & (52.6,42.4) \\
\hline
\end{tabular}
\caption{(ACC, NMI) Means by $L_{2,1}$ SNF with $(\alpha, \beta)$ on USPST}
\label{ab_uspst}
\end{table}

\subsection{Convergence Analysis for $L_{2,1}$ SNF}
\label{subsec:convergence}

The updating rules for minimizing the objective function of $L_{2,1}$ SNF are essentially iterative. We have shown that these update rules yield convergence, and we now investigate the rapidity of this convergence. To run the experiments, we set $k=5$ for Ionosphere, $k=10$ for Waveform, and $k=16$ for USPST, respectively. Meanwhile, the optimal $(\alpha, \beta)$ pair identified in Subsection \ref{subsec:parameter} for each dataset is adopted.

\begin{figure}[htb]
\hspace{-10mm}
\subfloat[Ionosphere]{\includegraphics[width=0.39\textwidth]{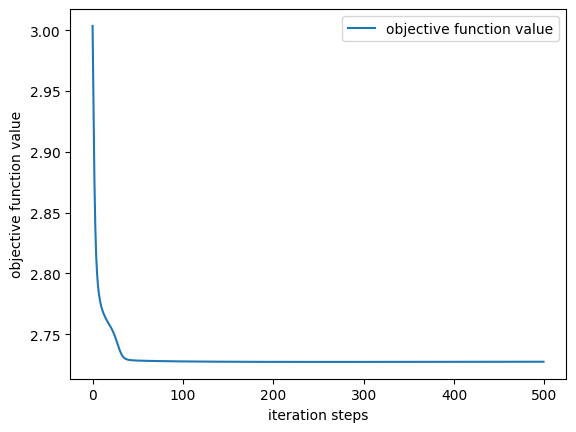}}
\subfloat[Waveform]{\includegraphics[width=0.39\textwidth]{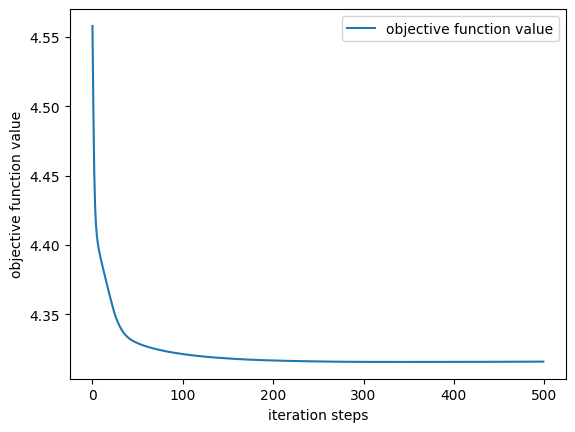}}
\subfloat[USPST]{\includegraphics[width=0.39\textwidth]{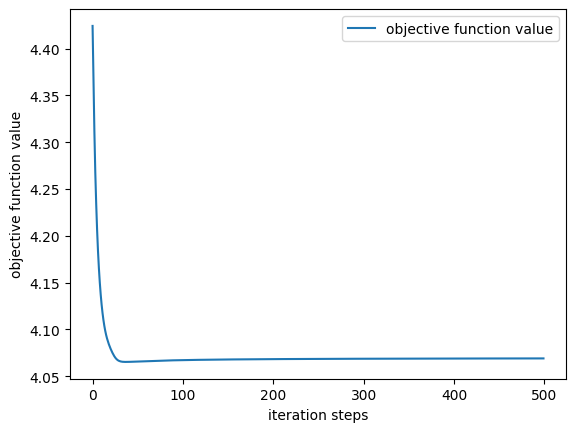}}
\caption{Convergence analysis of $L_{2,1}$ SNF on: (a) Ionosphere; (b) Waveform; and (c) USPST.
The y-axis for the objective function value is in the log scale.}
\label{convergencefig}
\end{figure}

Figure~\ref{convergencefig} shows the convergence curves of $L_{2,1}$ SNF on all three datasets with the optimal $\alpha$ and $\beta$. For each dataset, $x$-axis denotes the iteration number and $y$-axis stands for the log value of the objective function. We can see that the multiplicative update rules converge very quickly, usually within dozens of iterations.

It is imperative to provide some clarification towards this experiment. Firstly, due to the small cluster size $k$ compared with the actual data size $m$ and $n$, the leading term involving matrix decomposition error and thus the whole objective function $\mathscr{J}(\mathbf{U},\mathbf{V})$ in \eqref{l21snf} are rather large. Secondly, since other SNF algorithms use different objective function formulas, it would be unrealistic to compare them together here. However, other algorithms also demonstrate similar rapid convergence \cite{Ding2010, Luo2017}. Thirdly, as there are no reported algorithms to optimize $\mathscr{J}(\mathbf{U},\mathbf{V})$ in general due to NP-hardness, the limit value found in this experiment cannot be regarded as the minimization value of \eqref{l21snf}. To further explore the relationship between the optimal value in \eqref{l21snf} and the limit obtained by $L_{2,1}$ SNF algorithm, we run another experiment in Subsection \ref{subsec:random} on randomized matrices where $\displaystyle{\min_{\mathbf{U},\mathbf{V}\geq 0}}\mathscr{J}(\mathbf{U},\mathbf{V})=0$ is enforced and we observe that $L_{2,1}$ SNF algorithm reaches the optimal value, no matter what $\mathbf{U}$, $\mathbf{V}$ and their initialization are chosen randomly in advance.

\subsection{Experiments on Randomized Mixed-Sign Matrices}
\label{subsec:random}

Based on the performance comparison of the parameters for $L_{2,1}$ SNF in Subsection \ref{subsec:parameter}, we observe that the leading term of matrix decomposition plays the most important role while the other two terms improve the algorithm performance by varying the parameters $\alpha$ and $\beta$. Therefore, for fair comparison of all SNF algorithms without the effect of tuning parameters, we set $\alpha=\beta=0$ for $L_{2,1}$ SNF and Graph Regularized SNF, which reduces to the original SNF algorithm, and compare the numerical performance of $L_{2,1}$ SNF algorithm versus SNF algorithm. To this end, we begin with a randomized mixed-sign matrix $\mathbf{X}$ of dimension $10,000 \times 128$ in the range $[-1,1]$. We perform different degrees of size reduction with $k=32$ and $k=16$, respectively. In fact, these extreme matrix dimensions are inspired by the potential applications towards highly over-determined systems from the text mining and genomic data compression. The initialization values of $\mathbf{U}$ and $\mathbf{V}$ are randomized in the range of $[-1,1]$ and $[0,1]$, respectively. We then compare the relative error $\frac{\|\mathbf{X}-\mathbf{UV}^T\|_{2,1}}{\|\mathbf{X}\|_{2,1}}$, which is more suitable than the absolute error $\|\mathbf{X-UV^T}\|_{2,1}$ without worrying about the range of entries of $\mathbf{X}$, for SNF and $L_{2,1}$ SNF algorithms at different iteration steps.

\begin{figure}[htb]
\hspace{-2mm}
\subfloat[$k=32$]{\includegraphics[width=0.55\textwidth]{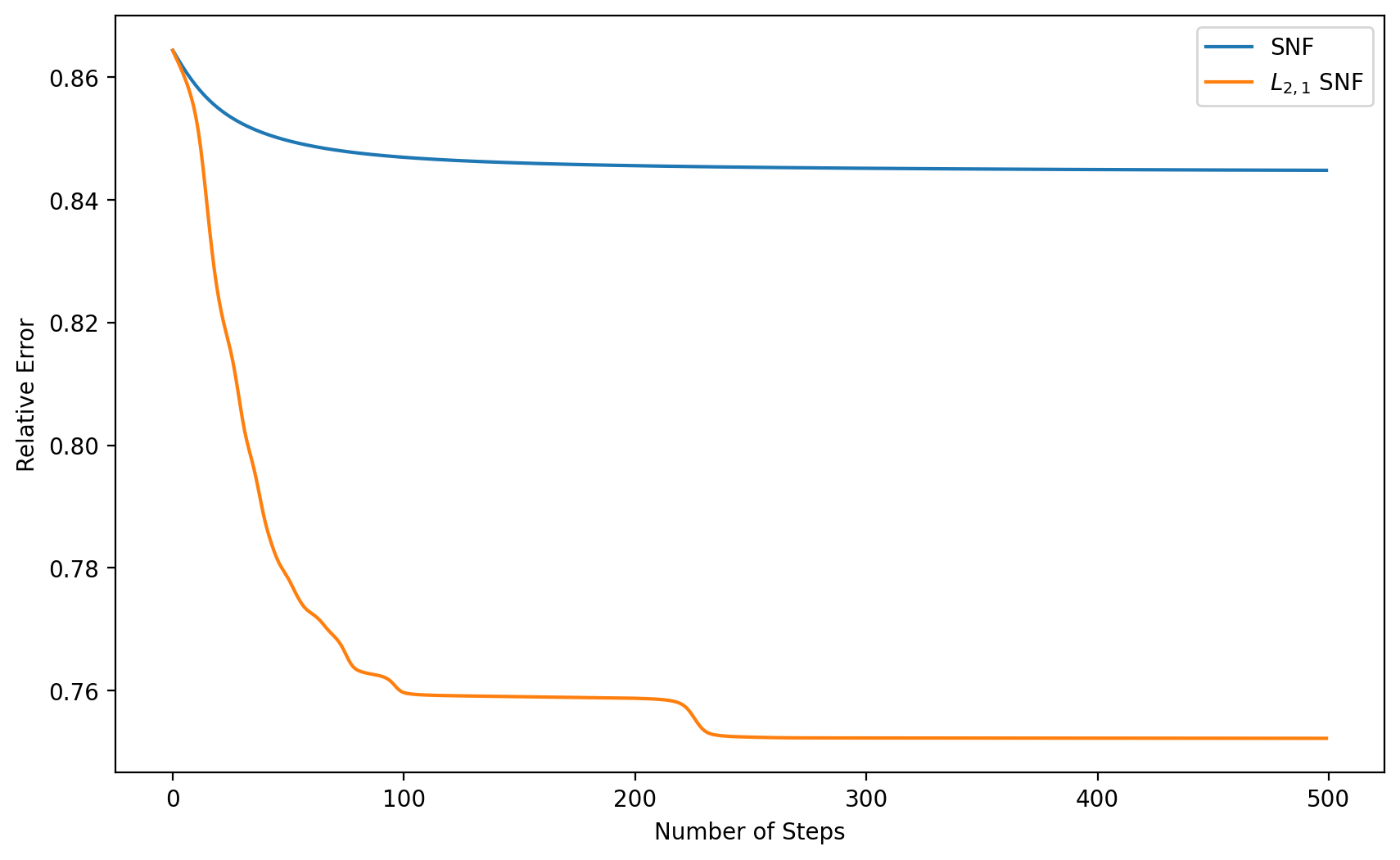}}
\subfloat[$k=16$]{\includegraphics[width=0.55\textwidth]{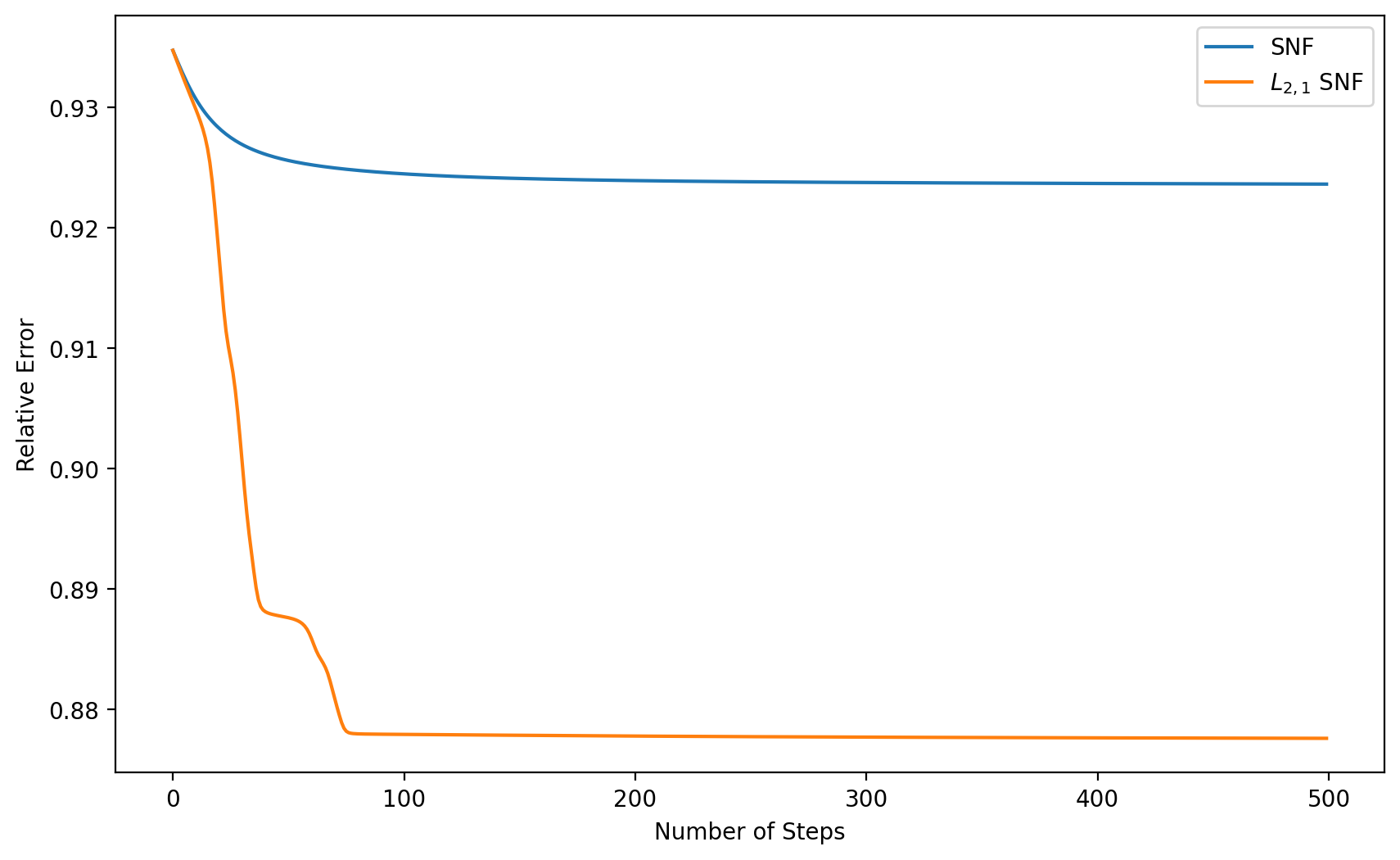}}
\caption{Comparison of the relative errors for $L_{2,1}$ SNF vs SNF algorithms for reduction of matrix $\mathbf{X}$ of dimension ${10,000 \times 128}$ with $k=32$ and $k=16$.}
\label{fig:matrixconvergence}
\end{figure}

Overall, $L_{2,1}$ SNF algorithm demonstrates a substantial improvement over SNF algorithm in reducing the relative errors, as shown in Figure \ref{fig:matrixconvergence}.

\begin{figure}[htb]
\hspace{-2mm}
\subfloat[$k=32$]{\includegraphics[width=0.55\textwidth]{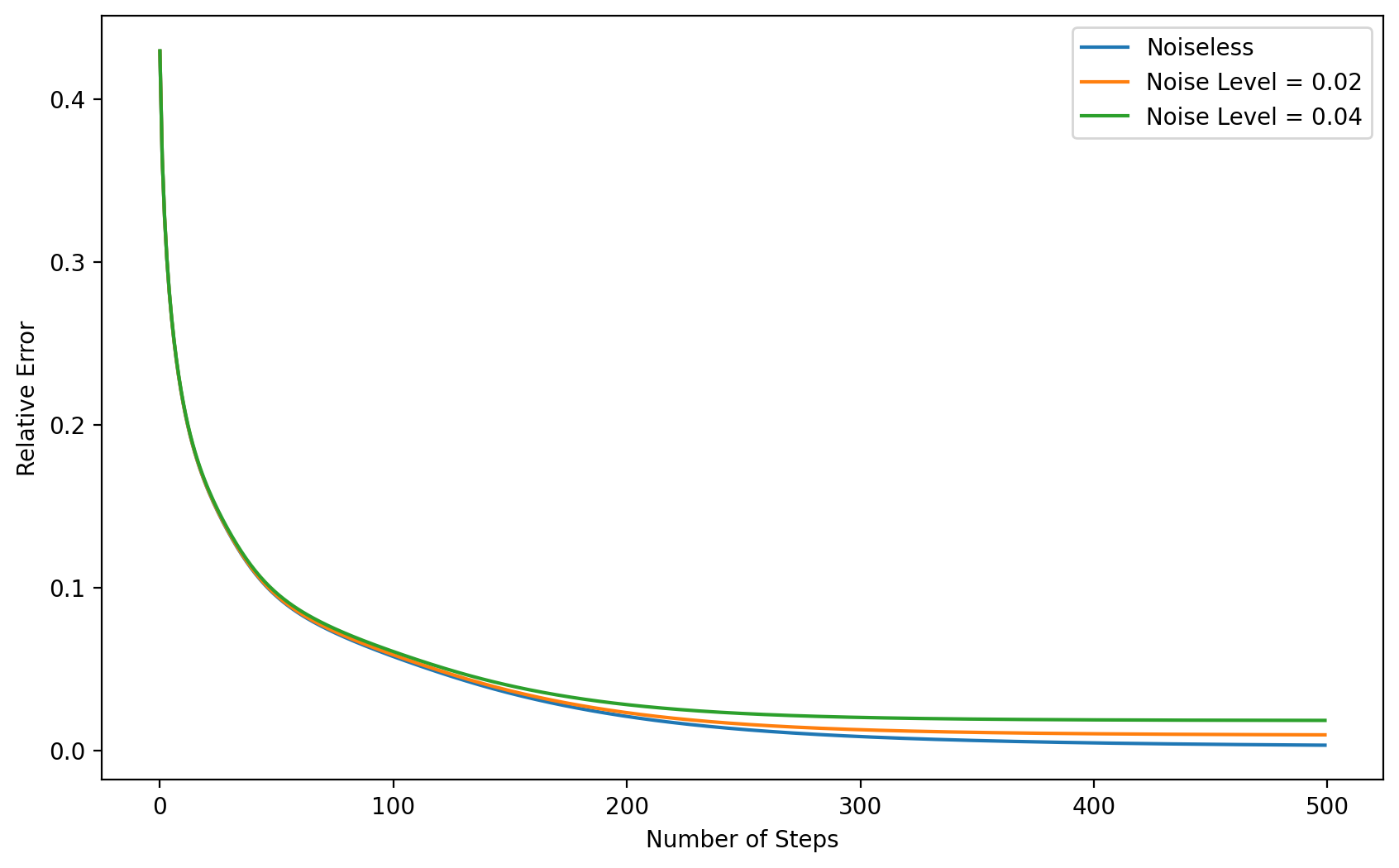}}
\subfloat[$k=16$]{\includegraphics[width=0.55\textwidth]{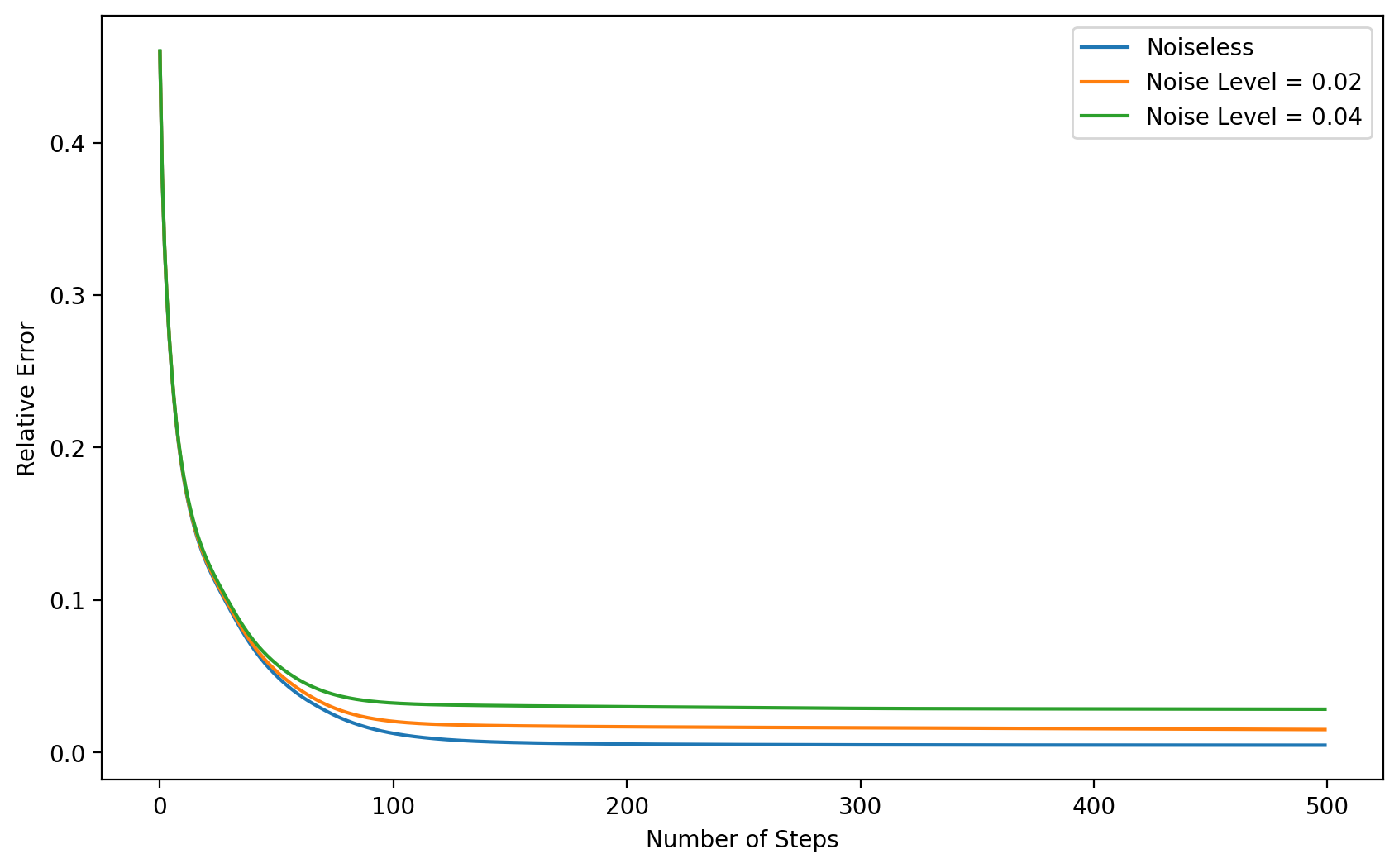}}
\caption{The relative errors of the object functions for $L_{2,1}$ SNF towards reduction of matrix $\mathbf{X}$ of dimension ${10,000 \times 128}$ with $k=32$ and $k=16$ at different levels of Gaussian noises.}
\label{fig:matrixtruesolution}
\end{figure}

Since it is unclear whether NMF and SNF algorithms converge to global minimizers of the objective functions theoretically, it is necessary to investigate how close their solutions are to the global minimizers in practice. Therefore, we run another experiment for $L_{2,1}$ SNF algorithm where $\mathbf{U}$ of dimension $10,000 \times k$ and $\mathbf{V}$ of dimension $128 \times k$ are generated randomly in the range of $[-1,1]$ and $[0,1]$, respectively for the cases of $k=32$ and $k=16$. We then enforce $\mathbf{X}=\mathbf{UV}^T$ of dimension $10,000 \times 128$ so that the objective function can indeed reach $\displaystyle{\min_{\mathbf{U}, \mathbf{V}\geq 0}}\|\mathbf{X} -\mathbf{UV}^T\|_{2,1}=0$ for all the cases. The initialization values for $\mathbf{U}$ and $\mathbf{V}$ are still randomized in the range of $[-1, 1]$ and $[0, 1]$, respectively. Furthermore, we add Gaussian noise at different levels with $\sigma=0.00$ (noiseless),  $\sigma=0.02$ and $\sigma=0.04$ to observe the robustness of $L_{2,1}$ SNF algorithm. We then compare the relative error $\frac{\|\mathbf{X}-\mathbf{UV}^T\|_{2,1}}{\|\mathbf{X}\|_{2,1}}$ at different iteration steps. The experimental results in Figure \ref{fig:matrixtruesolution} show that without Gaussian noise, the objective function always converges to the optimal value $0$, irrespective of the randomized initial values of $\mathbf{U}$ and $\mathbf{V}$. When the level of Gaussian noise increases, the error becomes larger but still stays within a reasonable range around $0$. In fact, original SNF algorithm also shows similar numerical behavior but we skip it here for brevity. This experiment confirms reported robustness of SNF algorithms. It also necessitates theoretical investigation on whether such limits from SNF algorithms can reach $\displaystyle{\min_{\mathbf{U}, \mathbf{V}\geq 0}}\|\mathbf{X}-\mathbf{UV}^T\|_{2,1}$.

\subsection{Experiments on CIFAR-10 Image Dataset}
\label{subsec:cifar10}

There exists many large-scale benchmark datasets used to measure the performance of new algorithms, such as CIFAR-10, Fashion-MNIST, etc. All such datasets turn out to contain nonnegative data which is the most popular data format in the field of machine learning. We now conduct a new experiment on the CIFAR-10 dataset which consists of 50000 training images and 10000 test images in 10 classes with corresponding labels in the range of $0\!-\!9$. Each image is a $32\times 32$ color image under the RGB color model and can thus be represented by a column vector of size $32\times 32\times 3 = 3072$ whose entries are in the range of $[0, 1]$. We construct a matrix $\mathbf{X}\geq 0$ of size $3072\times 10000$ to represent the whole set of test images by columns and then apply NMF algorithm \eqref{nmfsol}, SNF algorithm \eqref{snfsol}, and $L_{2,1}$ SNF algorithm with $\alpha=\beta=0$ \eqref{l21snfsol} towards decomposition of $\mathbf{X}$. We perform different degrees of size reduction with $k=15$ and $k=30$, respectively. The initialization values of $\mathbf{U}$ and $\mathbf{V}$ are randomized in the range of $[0,1]$. We compare the relative error $\frac{\|\mathbf{X}-\mathbf{UV}^T\|_{2,1}}{\|\mathbf{X}\|_{2,1}}$ for all the algorithms at different iteration steps.

As we know, SNF algorithms \eqref{snfsol} and \eqref{l21snfsol} are designed exclusively for mixed-sign data decomposition. Then for the nonnegative image matrix $\mathbf{X}$, SNF algorithms yield mixed-sign matrix $\mathbf{U}$ and may lead to undesired data reconstruction of the images due to existence of negative values, while NMF algorithm \eqref{nmfsol} works seamlessly due to its $\mathbf{U}\geq 0$ constraint. Therefore, it would be unrealistic to compare ACC and NMI based on the coordinate matrix $\mathbf{V}$. Nevertheless, it is still reasonable to compare their relative errors to observe which algorithm yields the best approximation of $\mathbf{X}$.

\begin{figure}[htb]
\hspace{-2mm}
\subfloat[$k=15$]{\includegraphics[width=0.55\textwidth]{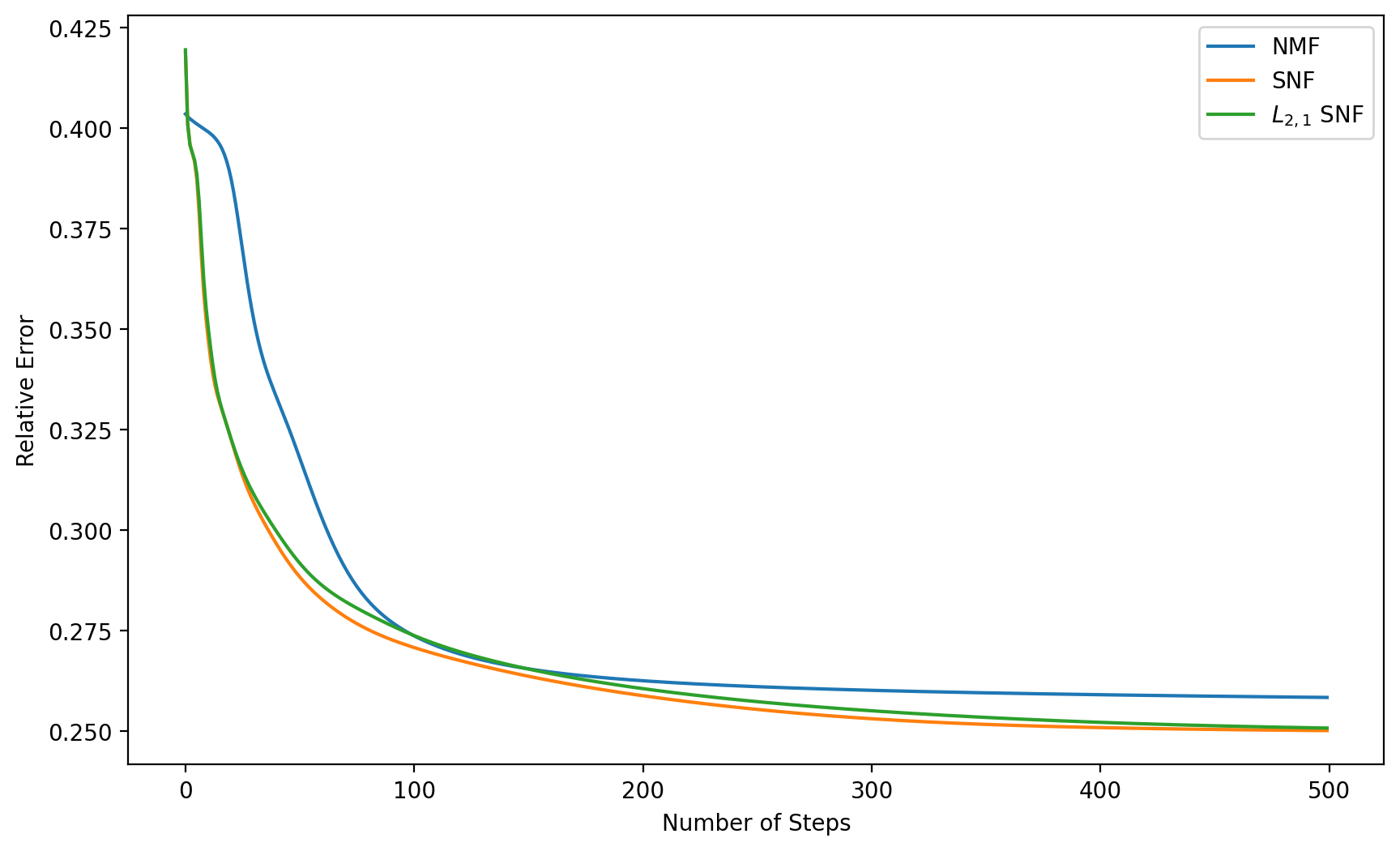}}
\subfloat[$k=30$]{\includegraphics[width=0.55\textwidth]{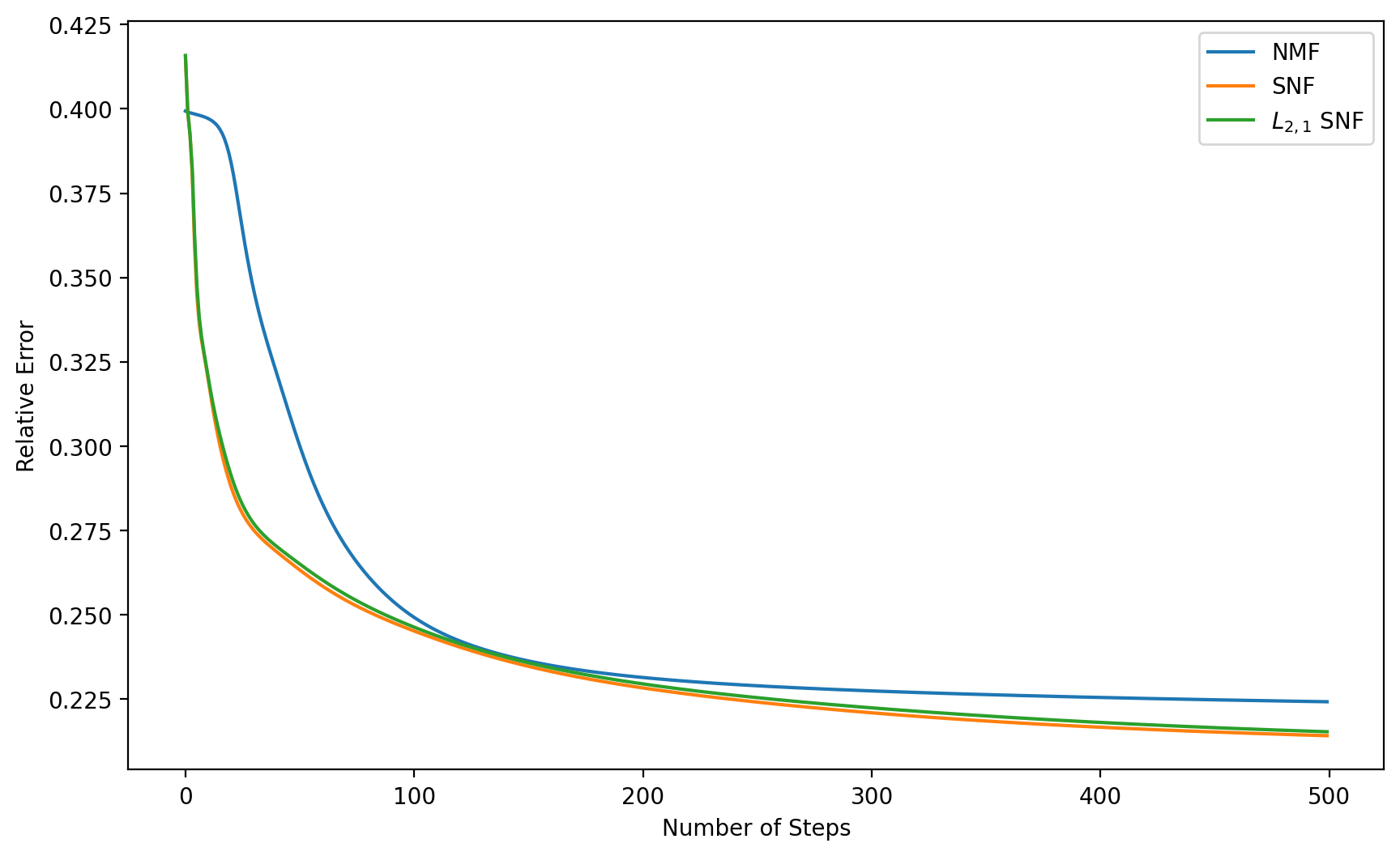}}
\caption{Comparison of the relative errors for NMF, SNF and $L_{2,1}$ SNF algorithms for reduction of CIFAR-10 test images of dimension ${3072 \times 10000}$ with $k=15$ and $k=30$.}
\label{fig:cifar10convergence}
\end{figure}

As shown in Figure \ref{fig:cifar10convergence}, SNF and $L_{2,1}$ SNF algorithms demonstrate better performance than NMF algorithm in reducing the relative errors. Meanwhile, $L_{2,1}$ SNF shows similar performance as original SNF as the proposed noise reduction mechanism is developed towards mixed-sign data. Therefore, this experiment verifies robustness of $L_{2,1}$ SNF algorithm \eqref{l21snfsol} in approximation of nonnegative datasets although it is designed for mixed-sign data.

\section{Conclusion and Future Research}
\label{sec:conclusion}

We have presented a novel data reduction algorithm, $L_{2,1}$ SNF, which renders a parts-based compression of mixed-sign data while reducing the effects of noise and outliers. We provide theoretical proof of monotonic convergence of the iterative updates given by the algorithm. Through experiments on three benchmark mixed-sign datasets and several randomized matrices with different levels of Gaussian noise, we demonstrate the advantage of $L_{2,1}$ SNF over other conventional SNF algorithms. In our future work, we will investigate theoretically the relationship between the solutions obtained by SNF algorithms and the global optimizers of the objective functions. We will also incorporate prior discriminative label information within $L_{2,1}$ SNF to construct a semi-supervised algorithm which is better-suited for real-world applications, including data mining and image processing.

\section*{Acknowledgement}
The authors thank the two anonymous reviewers and the editor for their useful suggestions which significantly improve the quality of this paper. This work is partially supported by NSF Research Training Group Grant DMS-2136228.

\section*{Conflict of Interest Statement}
The authors declare that they have no conflict of interest.

\bibliography{references}

\end{document}